\input kt3.mac

\input macro.mac


\null

\vskip 0.8 true cm

\centerline{\bf Nonlinear scalar field equations with $L^2$ constraint:}
\centerline{\bf Mountain pass and symmetric mountain pass approaches}

\bigskip

\centerline{Jun Hirata and Kazunaga Tanaka\footnote{${}^*$}{The second author is partially supported by JSPS Grants-in-Aid for Scientific 
Research (B) (25287025) and (B) (17H02855).}}

\medskip

\settabs 25 \columns
\+&& Department of Mathematics, School of Science and Engineering\cr
\+&& Waseda University, 3-4-1 Ohkubo, Shinjuku-ku, Tokyo 169-8555, Japan\cr

\bigskip


{\narrower

\noindent
{\bf Abstract:} 
We study the existence of radially symmetric solutions of the following nonlinear scalar field 
equations in $\R^N$ ($N\geq 2$):
    $$  (*)_m \left\{
        \eqalign{
        -&\Delta u = g(u) -\mu u \quad \hbox{in}\ \R^N, \cr
        &\norm u_{L^2(\R^N)} = m, \cr
        &u \in H^1(\R^N), \cr} \right.
    $$
where $g(\xi)\in C(\R,\R)$, $m>0$ is a given constant and $\mu\in \R$ is a Lagrange multiplier.

We introduce a new approach using a Lagrange formulation of the problem $(*)_m$.
We develop a new deformation argument under a new version of the Palais-Smale condition.
For a general class of nonlinearities related to [\cite[BL1], \cite[BL2], \cite[HIT]],
it enables us to apply minimax argument for $L^2$ constraint problems and we show the existence
of infinitely many solutions as well as mountain pass characterization of a minimizing solution
of the problem:
    $$  \inf\{ \int_{\R^N} \half \abs{\nabla u}^2 - G(u)\, dx;\, \norm u_{L^2(\R^N)}^2 = m \},
        \quad G(\xi)=\int_0^\xi g(\tau)\, d\tau.
    $$

}

\medskip

\BS{\label[Section:0]. Introduction}
In this paper, we study the existence of radially symmetric solutions of the following nonlinear
scalar field equations in $\R^N$ ($N\geq 2$):
    $$  (*)_m\ \left\{
        \eqalign{
        -&\Delta u = g(u) -\mu u \quad \hbox{in}\ \R^N, \cr
        &\norm u_{L^2(\R^N)}^2 = m, \cr
        &u \in H^1(\R^N), \cr} \right.
    $$
where $g(\xi)\in C(\R,\R)$, $m>0$ is a given constant and $\mu\in \R$ is a Lagrange multiplier.

Solutions of $(*)_m$ can be characterized as critical points of the constraint problem:
    $$  \calF(u)= \half\intRN \abs{\nabla u}^2 -\intRN G(u) :\, S_m \to \R,
    $$
where $S_m =\{ u\in \E;\, \norm u_{L^2(\R^N)}^2=m\}$ and $G(\xi)=\int_0^\xi g(\tau)\, d\tau$.

When $g(\xi)$ has $L^2$-subcritical growth, Cazenave-Lions [\cite[CL]] and Shibata [\cite[S1]] successfully 
found a solution of $(*)_m$ via minimizing method:
    $$  \calI_m = \inf_{u\in S_m} \calF(u).     \eqno\label[0.1]
    $$
[\cite[CL]] dealt with $g(\xi)=\abs\xi^{q-1}\xi$
($1<q<1+{4\over N}$) and [\cite[S1]] dealt with a class of more general nonlinearities, which satisfy conditions:

\itemitem{\newcond[g:1]} $g(\xi)\in C(\R,\R)$,
\itemitem{\newcond[g:2]} $\lim_{\xi\to 0} {g(\xi)\over \xi}= 0$,
\itemitem{\newcond[g:3]} $\lim_{\abs\xi\to\infty} {\abs{g(\xi)}\over \abs\xi^p}= 0$, where $p=1+{4\over N}$,
\itemitem{\newcond[g:4]} There exists $\xi_0>0$ such that $G(\xi_0)>0$.

\smallskip

\noindent
[\cite[S1]] showed

\smallskip

\item{(i)} There exists $m_S\geq 0$ such that for $m>m_S$.  $\calI_m$ defined in \ref[0.1] is achieved and
$(*)_m$ has at least one solution for $m>m_S$.
\item{(ii)} $m_S=0$ if and only if 
    $$  \lim_{\xi\to 0} {g(\xi)\over \abs\xi^{4\over N}\xi} = \infty.       \eqno\label[0.2]
    $$

\smallskip

\noindent
We remark that in [\cite[S1], \cite[CL]] they also studied orbital stability of the minimizer.
We also refer to Jeanjean [\cite[J]] and Bartsch-de Valeriola [\cite[BV]] for the study of
$L^2$-supercritical case (e.g. $g(\xi)\sim \abs{\xi}^{p-1}\xi$ with
$p\in (1+{4\over N},{N+2\over N-2})$).

We note that the conditions \cond[g:1]--\cond[g:4] are related to those in Berestycki-Lions [\cite[BL1], \cite[BL2]] 
(see also [\cite[BGK], \cite[HIT]]) as {\it almost necessary and sufficient conditions} for the existence 
of solutions of nonlinear scalar field equations:
    $$  \left\{ \eqalign{
        -&\Delta u = g(u) \quad \hbox{in}\ \R^N, \cr
        &u\in H^1(\R^N).  \cr}
        \right.                 \eqno\label[0.3]
    $$
More precisely, replacing \cond[g:2] by $\limsup_{\xi\to 0} {g(\xi)\over \xi}<0$ and replacing $p$ 
with ${N+2\over N-2}$ in \cond[g:3], they showed the existence of a least energy solution and 
they also showed the existence of a unbounded sequence of possibly sign-changing solutions 
assuming oddness of $g(\xi)$ in addition:

\smallskip

\itemitem{\newcond[g:5]}  $g(-\xi)=-g(\xi)$ for all $\xi\in \R$.

\smallskip

\noindent
We remark that if $g(\xi)$ satisfies \cond[g:1]--\cond[g:4], then $\tg(\xi)=g(\xi)-\mu\xi$ satisfies
the conditions of [\cite[BL1], \cite[BL2]] for $\mu\in (0,\infty)$ small.

In [\cite[CL], \cite[S1]], to show the achievement of $\calI_m$ on $S_m$, the following {\it sub-additivity inequality}
plays an important role.
    $$  \calI_m < \calI_s + \calI_{m-s} \quad \hbox{for all}\ s\in (0,m),               \eqno\label[0.4]
    $$
which ensures compactness of minimizing sequences for $\calI_m$.

In this paper, we take another approach to $(*)_m$ and we try to apply minimax methods to a Lagrange
formulation of the problem $(*)_m$:
    $$  
    \calL(\mu,u)=\half\intRN \abs{\nabla u}^2 -\intRN G(u) +{\mu\over 2}\left(\norm u_{L^2(\R^N)}^2 -m\right)
            :\, \RE\to \R.          \eqno\label[0.5]
    $$
We give another proof to the existence result of [\cite[S1]]; we take an approach related to 
Hirata-Ikoma-Tanaka [\cite[HIT]] and Jeanjean [\cite[J]], which made use of the scaling properties
of the problems to generate Palais-Smale sequences in augmented spaces with extra properties related to 
the Pohozaev identities.   We remark that such approaches were successfully applied to other problems 
with suitable scaling properties.  See 
Azzollini-d'Avenia-Pomponio [\cite[AdAP]], Byeon-Tanaka [\cite[BT]], Chen-Tanaka [\cite[CT]]
and Moroz-Van Schaftingen [\cite[MVS]].
We also give a mountain pass characterization of the minimizing value $\calI_m$ through the 
functional \ref[0.5], which we expect to be useful in the study of singular perturbation problems.  
We remark that a mountain pass characterization of the least energy solutions
for nonlinear scalar field equations \ref[0.3] was given in Jeanjean-Tanaka [\cite[JT]].

\medskip

\proclaim Theorem \label[Theorem:0.1].  
Assume \cond[g:1]--\cond[g:4].  Then we have
\item{(i)} There exists $m_0\in [0,\infty)$ such that for $m>m_0$, $(*)_m$ has at least one solution.
\item{(ii)} Assume \ref[0.2] in addition to  \cond[g:1]--\cond[g:4].  Then $(*)_m$ has at least one solution
for all $m>0$. 
\item{(iii)} In the setting of (i)--(ii),  a solution is obtained through a mountain pass minimax method: 
    $$  b_{mp} = \inf_{\gamma\in\Gamma_{mp}} \max_{t\in [0,1]} I(\gamma(t)).
    $$
See Section \ref[Section:5] for a precise definition of the minimax class $\Gamma_{mp}$.  We also have
    $$  b_{mp} = \calI_m,
    $$
where $\calI_m$ defined in \ref[0.1].

\medskip

\noindent
We will give a presentation of $m_0$ using least energy levels of $-\Delta u+\mu u =g(u)$ in Section \ref[Section:5].
We also show $m>m_0$ if and only if $\calI_m<0$.

We also deal with the existence of infinitely many solutions assuming oddness of $g(\xi)$.  It seems that
the existence of infinitely many solutions for the $L^2$-constraint problem is not well-studied.  
Our main result is the following

\medskip

\proclaim Theorem \label[Theorem:0.2].
Assume \cond[g:1]--\cond[g:4] and \cond[g:5].  Then
\item{(i)} For any $k\in\N$ there exists $m_k\geq 0$ such that for $m>m_k$, $(*)_m$ has at least $k$ solutions.
\item{(ii)} Assume \ref[0.2] in addition to \cond[g:1]--\cond[g:5].  Then for any $m>0$, $(*)_m$ has 
countably many solutions $(u_n)_{n=1}^\infty$, which satisfies
    $$  \eqalign{
        &\calF(u_n) < 0 \quad \hbox{for all}\ n\in\N, \cr
        &\calF(u_n) \to 0 \quad \hbox{as}\ n\to\infty. \cr}
    $$

\medskip

\noindent
To show Theorem \ref[Theorem:0.2], we develop a version of symmetric mountain pass methods, in which genus plays
an important role.

In the following sections, we give proofs to our Theorems \ref[Theorem:0.1] and \ref[Theorem:0.2].  Since the existence part of
Theorem \ref[Theorem:0.1] is already known by [\cite[CL], \cite[S1]].  We mainly deal with Theorem \ref[Theorem:0.2]
in Sections \ref[Section:1]--\ref[Section:4].

In Section \ref[Section:1], first we give a variational formulation of the problem $(*)_m$.  
For a technical reason, we write $\mu=e^\lambda$ ($\lambda\in\R$) and we try to find  critical points of 
    $$  I(\lambda,u)=\intRN \half\abs{\nabla u}^2 - G(u) +{e^\lambda\over 2}\left(\intRN\abs u^2-m\right)
        \in C^1(\RE,\R).
    $$
We also setup function spaces. Second for a fixed $\lambda\in\R$, we study the symmetric mountain pass value
$a_k(\lambda)$ of 
    $$  u\mapsto \whI(u) = \intRN \half\abs{\nabla u}^2 +{e^\lambda\over 2}u^2 - G(u).
    $$
Behavior of $a_k(\lambda)$ is important in our study.  In particular, $m_k$ in Theorem \ref[Theorem:0.2] is given by
    $$  m_k = 2 \inf_{\lambda\in (-\infty,\lambda_0)} {a_k(\lambda)\over e^\lambda}.
    $$
In Sections \ref[Section:2]--\ref[Section:3], we find that $I(\lambda,u):\,\RE\to \R$ has a kind of mountain pass geometry and we give a 
family of minimax sets for $I(\lambda,u)$, which involves the notion of genus under $\Z_2$-invariance:
$I(\lambda,-u)=I(\lambda,u)$.

In Section \ref[Section:4], we develop a new deformation argument to justify the minimax methods in Section \ref[Section:2].
Usually deformation theories are developed under so-called Palais-Smale condition.  However, under the conditions
\cond[g:1]--\cond[g:4], it is difficult to check the standard Palais-Smale condition for $I(\lambda,u)$.
We introduce a new version of Palais-Smale condition, which is inspired by our earlier work [\cite[HIT]] and
Jeanjean [\cite[J]].  Here we extend this idea further to generate a deformation flow, which is different from
the standard one;  our flow does not come from ODE in $\RE$ and in general it is not of class $C^1$.

In Section \ref[Section:5], we deal with Theorem \ref[Theorem:0.1] and we study the minimizing problem \ref[0.1].  
Applying mountain pass approach to $I(\lambda,u)$, we give another proof of the existence result 
as well as a mountain pass characterization of $\calI_m$ in $\RE$.

Our new deformation argument also works for the nonlinear scalar field equations \ref[0.3]; we can give a 
simpler proof to the results in [\cite[HIT]].  We believe that it is of interest and the idea is applicable 
to other problems with scaling properties.

\medskip

\BS{\label[Section:1]. Preliminaries}
\BSS{\label[Subsection:1.a].  Functional settings}
In Sections \ref[Section:1]--\ref[Section:4], we deal with Theorem \ref[Theorem:0.2] and we assume \cond[g:1]--\cond[g:5].
We denote by $\E$ the space of radially symmetric functions $u(x)=u(\abs x)$ which satisfy
$u(x)$, $\nabla u(x)\in L^2(\R^N)$.  We also use notation 
    $$  \eqalign{
        &\norm u_r = \left(\intRN \abs{u(x)}^r\right)^{1/r}  \quad \hbox{for}\ r\in [1,\infty) \ \hbox{and}\ u\in L^r(\R^N), \cr
        &\norm u_{H^1} = (\norm{\nabla u}_2^2 + \norm u_2^2)^{1/2} \quad \hbox{for}\ u\in \E. \cr}
    $$
We also write
    $$  (u,v)_2 = \intRN uv \quad \hbox{for}\ u,v\in L^2(\R^N).
    $$
In what follows, we denote by $p$ the $L^2$ critical exponent, i.e.,
    $$  p = 1 +{4\over N}.      \eqno\label[1.1]
    $$
In particular we have
    $$  {p+1\over p-1}-{N\over 2} = 1,      \eqno\label[1.2]
    $$
which we will use repeatedly in this paper.

For technical reasons, we set $\mu=e^\lambda$ in \ref[0.5] and we set for a given
$m>0$
    $$  I(\lambda,u)=\half\norm{\nabla u}_2^2 -\int_{\R^N} G(u) +{e^\lambda\over 2}\left(\norm u_2^2 -m\right) 
            :\, \RE\to \R.
    $$
It is easy to see that $I(\lambda,u)\in C^1(\RE,\R)$ and solutions of $(*)_m$
can be characterized as critical points of $I(\lambda,u)$, that is, $(\mu,u)$ with $\mu=e^\lambda>0$ solves $(*)_m$ 
if and only if $\partial_\lambda I(\lambda,u)=0$ and $\partial_u I(\lambda,u)=0$.
We also have
    $$  I(\lambda,-u) = I(\lambda,u) \quad \hbox{for all}\ (\lambda,u)\in \RE.    \eqno\label[1.3]
    $$
The following functionals will play important roles in our argument.
    $$  \eqalignno{
        \whI(\lambda,u) &=\half\norm{\nabla u}_2^2 + {e^\lambda\over 2}\norm u_2^2 -\intRN G(u):\, 
        \RE\to \R,                  &\label[1.4]\cr
        P(\lambda,u) &={N-2\over 2}\norm{\nabla u}_2^2 +N\left( {e^\lambda\over 2}\norm u_2^2 -\intRN G(u)\right):\,
        \RE\to \R.                  &\label[1.5]\cr}
    $$
We note that

\item{(i)} For a fixed $\lambda\in \R$, $u\mapsto \whI (\lambda,u)$ is corresponding to 
    $$  \left\{ \eqalign{
        -&\Delta u + e^\lambda u = g(u) \quad \hbox{in}\ \R^N, \cr
        &u\in \E. \cr}
        \right.         \eqno\label[1.6]
    $$
It is easy to see that
    $$  I(\lambda,u)=\whI(\lambda,u) -{e^\lambda\over 2}m \quad \hbox{for all}\ (\lambda,u).    \eqno\label[1.7]
    $$
\item{(ii)} $P(\lambda,u)$ is related to the Pohozaev identity.  It is well-known that for $\lambda\in \R$ 
if $u(x)\in \E$ solves \ref[1.6], then $P(\lambda,u)=0$.

\medskip

\BSScap{\label[Subsection:1.b]. Some estimates for $\whI(\lambda,u)$}{\raw[1.b]. Some estimates for I}
First we observe that for $\lambda\ll 0$, $u\mapsto \whI(\lambda,u)$ satisfies the assumptions of
[\cite[BL1], \cite[BL2], \cite[BGK], \cite[HIT]] and possesses the symmetric mountain pass geometry.  In what follows,
we write
    $$  S^{k-1} = \{\xi\in\R^k;\, \abs\xi = 1\},\quad D_k = \{\xi\in\R^k;\, \abs\xi \leq 1\}.
    $$
and set
    $$  \lambda_0 = \cases{
        \log\left(\displaystyle 2\sup_{\xi\not=0} {G(\xi)\over\xi^2}\right) &if $\displaystyle\sup_{\xi\not=0} {G(\xi)\over\xi^2} <\infty$, \cr
        \infty &if $\displaystyle \sup_{\xi\not=0} {G(\xi)\over\xi^2} =\infty$. \cr}        \eqno\label[1.8]
    $$

We have

\medskip

\proclaim Lemma \label[Lemma:1.1].
\item{(i)} For $\lambda\in (-\infty,\lambda_0)$,
    $$  G(\xi_0) -{e^\lambda\over 2}\abs{\xi_0}^2 > 0 \quad \hbox{for some}\ \xi_0>0.
    $$
In particular, $\widetilde g(\xi)=g(\xi)-e^\lambda \xi$ satisfies the assumptions of [\cite[BL1],\cite[BL2],\cite[HIT]],
that is, $\widetilde g(\xi)$ satisfies \cond[g:1], \cond[g:3]--\cond[g:5] and 
    $$  \lim_{\xi\to 0} {\tg(\xi)\over \xi}<0.
    $$
\item{(ii)} For any $\lambda\in (-\infty,\lambda_0)$ and for any $k\in \N$, there exists a continuous odd map 
$\zeta:\, S^{k-1}\to \E$ such that
    $$  \whI(\lambda,\zeta(\xi)) < 0 \quad \hbox{for all}\ \xi\in S^{k-1}.
    $$
\item{(iii)} When $\lambda_0<\infty$, for $\lambda\geq\lambda_0$ we have
    $$  G(\xi) -{e^\lambda\over 2}\abs{\xi}^2 \leq 0 \quad \hbox{for all}\ \xi\in\R.
    $$
In particular, $\whI(\lambda,u)\geq 0$ for all $u\in\E$.

\medskip

\claim Proof.
By \cond[g:1]--\cond[g:5] and the definition \ref[1.8] of $\lambda_0$, we can easily see (i) and (iii). 
By the arguments in [\cite[BL2], \cite[HIT]], we can observe that $u\mapsto \whI(\lambda,u)$ has
the property (ii).  \QED

\medskip

\noindent
For $k\in\N$ and $\lambda\in (-\infty,\lambda_0)$, we set
    $$  \eqalignno{
        \whGamma_k(\lambda) &= \{ \zeta\in C(D_k,\E);\, \zeta(-\xi)=-\zeta(\xi) \ \hbox{for}\ \xi\in D_k, \cr
            &\qquad \whI(\lambda,\zeta(\xi)) < 0 \ \hbox{for}\ \xi\in \partial D_k=S^{k-1}\},   &\label[1.9]\cr
        \a_k(\lambda) &= \inf_{\zeta\in \whGamma_k(\lambda)} \max_{\xi\in D_k} \whI(\lambda,\zeta(\xi)).
                                                                                        &\label[1.10]\cr}
    $$
We note that $\whGamma_k(\lambda)\not=\emptyset$ by Lemma \ref[Lemma:1.1] (ii).  Since 
$\whI(\lambda,u) =\half (\norm{\nabla u}_2^2+e^\lambda\norm u_2^2) + o(\norm u_{H^1}^2)$ as $\norm u_{H^1}\sim 0$,
    $$  \a_k(\lambda) >0 \quad \hbox{for all}\ \lambda\in (-\infty,\lambda_0) \ \hbox{and}\ k\in \N.
    $$
By the results of [\cite[HIT]], we observe that $a_k(\lambda)$ is a critical value of $u\mapsto \whI(\lambda,u)$.
See also Section \ref[Section:6].

We also have
    $$  \eqalignno{
    &0<\a_1(\lambda) \leq \a_2(\lambda) \leq \cdots \leq\a_k(\lambda) \leq \a_{k+1}(\lambda) \leq \cdots 
        \quad \hbox{for all}\ \lambda\in (-\infty,\lambda_0),   &\label[1.11]\cr
    &\a_k(\lambda) \leq \a_k(\lambda') \quad \hbox{for all}\ \lambda<\lambda' < \lambda_0 \ \hbox{and} \ k\in \N.\cr}
    $$
For the behavior of $a_k(\lambda)$ as $\lambda\to -\infty$, the condition \ref[0.2] is important.  
We have

\medskip

\proclaim Lemma \label[Lemma:1.2].
Assume \cond[g:1]--\cond[g:5].
\item{(i)} Assume \ref[0.2] in addition.   Then for any $k\in \N$
    $$  \lim_{\lambda\to -\infty} {\a_k(\lambda)\over e^\lambda}=0.
    $$
\item{(ii)}  If 
    $$  \limsup_{\xi\to 0}{\abs{g(\xi)}\over \abs\xi^p} <\infty,        \eqno\label[1.12]
    $$
then for any $k\in\N$
    $$  \liminf_{\lambda\to -\infty} {a_k(\lambda)\over e^\lambda} > 0.
    $$

\claim Proof.
(i) Choose $r\in(1+{4\over N},{N+2\over N-2})$.  By \ref[0.2] and \cond[g:3], for any $L>0$ there exists $C_L>0$
such that
    $$  \xi g(\xi)\geq L\abs\xi^{p+1}-C_L\abs\xi^{r+1} \quad \hbox{for all}\ \xi\in \R,
    $$
from which we have
    $$  \eqalign{
        G(\xi) &\geq {L\over p+1}\abs\xi^{p+1} - {C_L\over r+1} \abs\xi^{r+1} 
                \quad \hbox{for all}\ \xi\in\R,\cr
        \whI(\lambda,u) &\leq \half\norm{\nabla u}_2^2 +{e^\lambda\over 2}\norm u_2^2 - {L\over p+1}\norm u_{p+1}^{p+1}
                    +{C_L\over r+1}\norm u_{r+1}^{r+1} \quad \hbox{for all}\ u\in \E.\cr}
    $$
Setting $u(x)=e^{\lambda\over p-1} v(e^{\lambda/2}x)$, $v(x)\in \E$, we have from \ref[1.2]
    $$  \whI(\lambda,u) \leq e^\lambda\left(\half\norm v_{H^1}^2 -{L\over p+1}\norm v_{p+1}^{p+1} 
            +{C_L\over r+1}e^{{r-p\over p-1}\lambda}\norm v_{r+1}^{r+1}\right).
    $$
We note that
    $$  \oI(v) = \half\norm v_{H^1}^2 -{1\over p+1}\norm v_{p+1}^{p+1}:\, \E\to \R
    $$
has the symmetric mountain pass geometry and thus there exists an odd continuous map 
$\ozeta(\xi):\, D_k\to \E$ such that $\oI(\ozeta(\xi)) <0$ for all $\xi\in \partial D_k$.
By \ref[1.7], $\zeta_\lambda(\xi)=e^{\lambda\over p-1}\ozeta(\xi)(e^{\lambda/2}x)$ satisfies for $L\geq 1$
    $$  \whI(\lambda,\zeta_\lambda(\xi)) \leq e^\lambda\left( \oI(\ozeta(\xi)) 
        -{L-1\over p+1}\norm{\ozeta(\xi)}_{p+1}^{p+1} 
        + {C_L\over r+1} e^{{r-p\over p-1}\lambda}\norm{\ozeta(\xi)}_{r+1}^{r+1}\right).
    $$
Thus for $\lambda\ll 0$, we have $\zeta_\lambda(\xi)\in \whGamma_k(\lambda)$ and we have
    $$  \limsup_{\lambda\to -\infty} {\a_k(\lambda)\over e^\lambda} \leq \max_{\xi\in D_k}
        \left(\half \norm{\ozeta(\xi)}_{H^1}^2 -{L\over p+1}\norm{\ozeta(\xi)}_{p+1}^{p+1}\right).
    $$
Since $L\geq 1$ is arbitrary, we have the conclusion.  

(ii) By \ref[1.12] and \cond[g:3], there exist $C>0$ such that
    $$  G(\xi) \leq C\abs\xi^{p+1} \quad \hbox{for all}\ \xi\in\R.
    $$
Thus we have
    $$  \whI(\lambda,u)\geq \half\norm{\nabla u}_2^2 +{e^\lambda\over 2}\norm u_2^2 -{C\over p+1}\norm u_{p+1}^{p+1}.
    $$
As in (i),
    $$  \whI(\lambda,e^{\lambda\over p-1}u(e^{\lambda/2}x)) 
        \geq e^\lambda\left(\half\norm{u(x)}_{H^1}^2 -{C\over p+1}\norm{u(x)}_{p+1}^{p+1}\right),
    $$
from which we deduce that $a_k(\lambda)/e^\lambda$ is estimated from below by the mountain pass minimax
value for $u\mapsto \half\norm{u(x)}_{H^1}^2 -{C\over p+1}\norm{u(x)}_{p+1}^{p+1}$.  Thus (ii) holds.  \QED

\medskip

\noindent
We define for $k\in \N$
    $$  m_k = 2 \inf_{\lambda\in (-\infty,\lambda_0)} 
        {\a_k(\lambda)\over e^\lambda} \geq 0.  \eqno\label[1.13]
    $$
By \ref[1.11], we have
    $$  0\leq m_1\leq m_2 \leq \cdots \leq m_k \leq m_{k+1} \leq \cdots.     \eqno\label[1.14]
    $$
In what follows, we fix $m>m_k$ arbitrary and try to show that $I(\lambda,u)$ has at least $k$ pairs of
critical points.

As a corollary to Lemma \ref[Lemma:1.2], we have

\medskip

\proclaim Corollary \label[Corollary:1.3]. 
\item{(i)}  Under the condition \ref[0.2],
    $$  m_k = 0 \quad \hbox{for all}\ k\in\N.
    $$
\item{(ii)} Under the condition \ref[1.12],
    $$  m_k >0 \quad \hbox{for all}\ k\in\N.                \QED
    $$

\medskip

\BSS{\label[Subsection:1.c].  An estimate from below}
By \cond[g:2] and \cond[g:3], for any $\delta>0$ there exists $C_\delta>0$ such that
    $$  \xi g(\xi) \leq C_\delta\abs\xi^2 + \delta\abs\xi^{p+1} \quad \hbox{for all}\ \xi\in \R.
                                                \eqno\label[1.15]
    $$
Then we also have
    $$  G(\xi) \leq \half C_\delta\abs\xi^2 + {\delta\over p+1}\abs\xi^{p+1} \quad \hbox{for all}\ \xi\in \R.
    $$
Setting 
    $$  \uI(\lambda,u) = \half\norm{\nabla u}_2^2 +\half(e^\lambda-C_\delta)\norm u_2^2 
        -{\delta\over p+1}\norm u_{p+1}^{p+1},  
    $$
we have
    $$  \whI(\lambda,u) \geq \uI(\lambda,u) \quad \hbox{for all}\ (\lambda,u)\in \RE.
                                                                \eqno\label[1.16]
    $$
The functional $u\mapsto \uI(\lambda,u)$ has a typical mountain pass geometry if $e^\lambda>C_\delta+1$, which
enables us to give an estimate of $\whI(\lambda,u)$ from below.

In what follows, we denote by $E_0>0$ the least energy level for $-\Delta u+u=\abs u^{p-1}u$ in $\R^N$, that is,
    $$  E_0 = \inf\{ \half\norm{\nabla u}_2^2 +\half \norm u_2^2 - {1\over p+1}\norm u_{p+1}^{p+1};\,
        u\not=0,\, \norm{\nabla u}_2^2 +\norm u_2^2 =\norm u_{p+1}^{p+1} \}.
    $$

\proclaim Lemma \label[Lemma:1.4].
For $e^\lambda\geq C_\delta +1$,
    $$  \eqalignno{
        \whI(\lambda,u) &\geq \delta^{-{2\over p-1}}(e^\lambda-C_\delta)E_0 \quad \hbox{if}\ u\not=0,\, 
            \norm{\nabla u}_2^2 +(e^\lambda-C_\delta)\norm u_2^2 =\delta\norm u_{p+1}^{p+1}, &\label[1.17]\cr
        \whI(\lambda,u) &\geq 0  \quad \hbox{if}\ 
            \norm{\nabla u}_2^2 +(e^\lambda-C_\delta)\norm u_2^2 \geq \delta\norm u_{p+1}^{p+1}. &\label[1.18]\cr}
    $$

\claim Proof.
Let $\omega(x)$ be the least energy solution of $-\Delta u+u = \abs u^{p-1}u$.  Then it is easy to see that
    $$  u_{\lambda,\delta}(x)= \left({e^\lambda-C_\delta\over \delta}\right)^{1\over p-1}\omega((e^\lambda-C_\delta)^{1/2} x)
    $$
is a least energy solution of $-\Delta u +(e^\lambda-C_\delta)u=\delta\abs u^{p-1}u$ in $\R^N$.
Set $S_{\lambda,\delta}=\{ u\in \E\setminus\{ 0\};\, \norm{\nabla u}_2^2 +(e^\lambda-C_\delta)\norm u_2^2 =\delta \norm u_{p+1}^{p+1}\}$.
By \ref[1.2], it is easy to see that for $e^\lambda\geq C_\delta +1$
    $$  \uI(\lambda,u) \geq \delta^{-{2\over p-1}}(e^\lambda-C_\delta)E_0 \quad \hbox{for}\ u \in S_{\lambda,\delta}.
    $$
Thus we get \ref[1.17] from \ref[1.16].  Noting
    $$  \{ u\in \E;\, \norm{\nabla u}_2^2 +(e^\lambda-C_\delta)\norm u_2^2 \geq \delta \norm u_{p+1}^{p+1}\}
        = \{ tu;\, t\in [0,1],\, u\in S_{\lambda,\delta} \}
    $$
and that for $u\in S_{\lambda,\delta}$, $\uI(\lambda,tu)$ is increasing for $t\in (0,1)$, we have \ref[1.18].  \QED


\medskip

\BScap{\label[Section:2]. Minimax methods for $I(\lambda,u)$}{2. Minimax methods for I}
\BSS{\label[Subsection:2.a]. Symmetric mountain pass methods}
We fix $k\in\N$ and $m>0$ such that
    $$  m > m_k,                                \eqno\label[2.1]
    $$
where $m_k\geq 0$ is given in \ref[1.13].  We will show that $I(\lambda,u)$ has at least $k$
pairs of critical points.

We choose $\delta_m>0$ such that
    $$  \delta_m^{-{2\over p-1}} E_0 > {m\over 2}               \eqno\label[2.2]
    $$
and take $C_{\delta_m}>0$ so that \ref[1.15] holds.  For 
    $$  \lambda_m = \log(C_{\delta_m}+1),
    $$
we set
    $$  \Omega_m=\hbox{int}\{ (\lambda,u);\, \lambda\geq \lambda_m,\, 
            \norm{\nabla u}_2^2+(e^\lambda-C_{\delta_m})\norm u_2^2 \geq \delta_m\norm u_{p+1}^{p+1}\}.
                                                        \eqno\label[2.3]
    $$
We note that $\Omega_m=\bigcup_{\lambda\in [\lambda_m,\infty)}(\{\lambda\}\times D_\lambda)$ is a domain
whose section $D_\lambda\subset\E$ is a set surrounded by the Nehari manifold $\{ u\in\E\setminus\{0 \};\, 
\norm{\nabla u}_2^2 +(e^\lambda-C_{\delta_m})\norm u_2^2 = \delta_m\norm u_{p+1}^{p+1}\}$.
In particular $[\lambda_m,\infty)\times\{ 0\}\subset\Omega_m$.

Using Lemma \ref[Lemma:1.4], we have

\medskip

\proclaim Lemma \label[Lemma:2.1].
\item{(i)} $\displaystyle B_m \equiv \inf_{(\lambda,u)\in\partial\Omega_m} I(\lambda,u) > -\infty$.
\item{(ii)} $\whI(\lambda,u)\geq 0$ for $(\lambda,u)\in \Omega_m$.

\medskip

\claim Proof.
Note that $\partial\Omega_m=\calC_0\cup \calC_1$, where
    $$  \eqalign{
        \calC_0 &= \{ (\lambda,u)\in \R\times(\E\setminus\{ 0\});\, \lambda\geq \lambda_m, \cr
                &\qquad \qquad \norm{\nabla u}_2^2 +(e^\lambda -C_{\delta_m})\norm u_2^2 =\delta_m \norm u_{p+1}^{p+1} \},\cr
        \calC_1 &= \{(\lambda_m,u)\in \RE;\, 
            \norm{\nabla u}_2^2 +(e^{\lambda_m} -C_{\delta_m})\norm u_2^2 \geq \delta_m \norm u_{p+1}^{p+1} \},\cr }
    $$
By Lemma \ref[Lemma:1.4],
    $$  \eqalign{
    I(\lambda,u) &= \whI(\lambda,u)-{e^\lambda\over 2}m 
        \geq \delta_m^{-{2\over p-1}}(e^\lambda-C_{\delta_m})E_0-{e^\lambda\over 2}m 
        \quad \hbox{for}\ (\lambda,u)\in \calC_0, \cr
    I(\lambda,u) &= \whI(\lambda,u)-{e^\lambda\over 2}m \geq -{e^{\lambda_m}\over 2}m 
        \quad \hbox{for}\ (\lambda,u)\in \calC_1. \cr}
    $$
By our choice \ref[2.2] of $\delta_m$, we have 
$\inf_{(\lambda,u)\in\partial\Omega_m} I(\lambda,u) > -\infty$ and (i) holds.
(ii) is also clear.         \QED

\medskip

We introduce a family of minimax methods.  For $j\in\N$ we set
    $$  \eqalign{
        \Gamma_j &= \{\gamma(\xi)=(\varphi(\xi), \zeta(\xi))\in C(D_j,\RE);\cr
        &\qquad \hbox{$\gamma(\xi)$ satisfies conditions \cond[\gamma:1]--\cond[\gamma:3] below}\},\cr}
    $$
where

\smallskip 

\itemitem{\newcond[\gamma:1]} $\varphi(-\xi)=\varphi(\xi)$, $\zeta(-\xi)=-\zeta(\xi)$ for all $\xi\in D_j$.
\itemitem{\newcond[\gamma:2]} There exists $\lambda\in (-\infty,\lambda_0)$ such that
    $$  \varphi(\xi)=\lambda, \ I(\lambda,\zeta(\xi)) \in (\RE)\setminus \Omega_m, \
            I(\lambda,\zeta(\xi)) \leq B_m-1 \quad \hbox{for}\ \xi\in \partial D_j.
    $$
\itemitem{\newcond[\gamma:3]} $\varphi(0) \in [\lambda_m,\infty)$ and $\zeta(0)=0$.  Moreover
    $$  I(\varphi(0), \zeta(0)) = -{e^{\varphi(0)}\over 2}m \leq B_m -1.
    $$

\smallskip

\noindent
We note that
\item{(i)} for $\lambda\in (-\infty,\lambda_0)$, $u\mapsto \whI(\lambda,u)$ has the symmetric mountain pass
geometry.
\item{(ii)} $I(\lambda,0)=-{e^\lambda\over 2}m \to -\infty$ as $\lambda\to \infty$.

\noindent
From these facts, we have $\Gamma_j\not=\emptyset$ for all $j\in\N$.

We remark that $\Gamma_j$ is a family of $j$-dimensional symmetric mountain paths joining points in
$[\lambda_m,\infty)\times\{ 0\}\subset\Omega_m$ and $(\RE)\setminus\Omega_m$.

We set
    $$  b_j = \inf_{\gamma\in\Gamma_j}\max_{\xi\in D_j} I(\gamma(\xi)) \quad \hbox{for}\ j\in\N.
    $$

\proclaim Proposition \label[Proposition:2.2].
\item{(i)} $b_j\geq B_m$ for all $j\in\N$.
\item{(ii)} $b_j<0$ for $j=1,2,\cdots, k$.

\medskip

To show Proposition \ref[Proposition:2.2], we need

\medskip

\proclaim Lemma \label[Lemma:2.3].
    $$  b_j \leq a_j(\lambda) -{e^\lambda\over 2}m \quad \hbox{for}\ \lambda\in (-\infty,\lambda_0).
                                                        \eqno\label[2.4]
    $$

\claim Proof.
First we note that by (ii) of Lemma \ref[Lemma:2.1] that
    $$  (\lambda,\zeta(\xi)) \in (\RE)\setminus\Omega_m \quad 
        \hbox{for}\ \zeta\in \whGamma(\lambda) \ \hbox{and}\ \xi\in \partial D_j.
    $$
Second we remark that we may assume for $\zeta(\xi)\in \whGamma_j(\lambda)$
    $$  I(\lambda,\zeta(\xi)) \leq B_m-1 \quad \hbox{for}\ \xi\in\partial D_j.   \eqno\label[2.5]
    $$
In fact, for $u\in \E$ and $\nu>0$ we have
    $$  \whI(\lambda,u(x/\nu)) = \half\nu^{N-2} \norm{\nabla u}_2^2
        +\nu^N\left({e^\lambda\over 2}\norm u_2^2 -\intRN G(u)\right),
    $$
from which we deduce that if $\whI(\lambda,u(x))<0$, then $\nu\mapsto \whI(\lambda,u(x/\nu));\, [1,\infty)\to \R$
is decreasing and $\lim_{\nu\to \infty}\whI(\lambda,u(x/\nu))=-\infty$. 
Thus, for a given $\zeta(\xi)\in \whGamma_j(\lambda)$, setting
    $$  \widetilde\zeta(\xi)(x) = \cases{
        \zeta(2\xi) &for $\abs\xi\in [0, \half]$, \cr
        \zeta(\xi/\abs\xi)({x\over L(2\abs\xi-1)+1}) &for $\abs\xi\in (\half,1]$.\cr}
    $$
We find for $L\gg 1$, $\widetilde \zeta(\xi)\in \whGamma_j(\lambda)$ and
    $$  \eqalign{
        &\max_{\xi\in D_j}\whI(\lambda,\widetilde\zeta(\xi)) =\max_{\xi\in D_j} \whI(\lambda,\zeta(\xi)),\cr
        &I(\lambda,\widetilde \zeta(\xi)) \leq B_m-1 \quad \hbox{for}\ \xi\in \partial D_j. \cr}
    $$
Thus we may assume \ref[2.5] for $\zeta(\xi)\in \whGamma_j(\lambda)$.

Next we show \ref[2.4].  For $\zeta(\xi)\in \whGamma_j(\lambda)$ with \ref[2.5], we set $\check\gamma(\xi)=(\check\varphi(\xi),
\check\zeta(\xi))$ by
    $$  \eqalign{
        &\check\varphi(\xi)=\cases{
            \lambda+R(1-2\abs\xi)   &for $\abs\xi\in [0,\half]$,\cr
            \lambda                 &for $\abs\xi\in (\half,1]$,\cr} \cr
        &\check\zeta(\xi)=\cases{
            0                       &for $\abs\xi\in [0,\half]$,\cr
            \zeta({\xi\over\abs\xi}(2\abs\xi-1))    &for $\abs\xi\in (\half,1]$.\cr} \cr}
    $$
Then for $R$ large, we have $\check\gamma(\xi)\in \Gamma_j$ and
    $$  \eqalign{
        I(\check\gamma(\xi)) &= I(\lambda+R(1-2\abs\xi),0) = -{e^{\lambda+R(1-2\abs\xi)}\over 2}m 
            \leq -{e^\lambda\over 2}m \quad \hbox{for}\ \abs\xi \in [0,\half],\cr
        I(\check\gamma(\xi)) &= I(\lambda,\check\zeta(\xi)) = \whI(\lambda,\check\zeta(\xi)) -{e^\lambda\over 2}m 
            \leq \max_{\xi\in D_j} \whI(\lambda,\zeta(\xi)) -{e^\lambda\over 2}m 
            \quad \hbox{for}\ \abs\xi\in (\half,1]. \cr }
    $$
Since $\zeta(\xi)\in \whGamma_j(\lambda)$ is arbitrary, we have \ref[2.4].  \QED

\medskip

Now we give a proof to Proposition \ref[Proposition:2.2].

\medskip

\claim Proof of Proposition \ref[Proposition:2.2].
(i) By \cond[\gamma:2], \cond[\gamma:3], we have
    $$  \gamma(\partial D_j)\cap \Omega_m =\emptyset \quad \hbox{and}\quad \gamma(0)\in \Omega_m
        \quad \hbox{for all}\ \gamma\in \Gamma_j.
    $$
Thus $\gamma(D_j)\cap\partial \Omega_m\not=\emptyset$ for all $\gamma\in\Gamma_j$ and it follows 
from Lemma \ref[Lemma:2.1] (i) that
    $$  \max_{\xi\in D_j} I(\gamma(\xi)) \geq \inf_{(\lambda,u)\in \partial\Omega_m} I(\lambda,u) \equiv B_m.
    $$
Since $\gamma\in \Gamma_j$ is arbitrary, we have (i).\m
(ii) By Lemma \ref[Lemma:2.3], for any $\lambda\in (-\infty,\lambda_0)$
    $$  {b_j\over e^\lambda} \leq {a_j(\lambda)\over e^\lambda} -{m\over 2}.
    $$
Since
    $$  2\inf_{\lambda\in (-\infty,\lambda_0)} \left({a_j(\lambda)\over e^\lambda}-{m\over 2}\right) 
        = m_j - m,
    $$
the conclusion (ii) follows from \ref[1.14] and \ref[2.1].  \QED

\medskip

In Section \ref[Section:3], we will see that $I(\lambda,u)$ satisfies a version of Palais-Smale type condition $(PSP)_b$
for $b<0$, which enables us to develop a deformation argument and to show $b_j$ ($j=1,2,\cdots,k$) are
critical values of $I(\lambda,u)$.  However to show multiplicity, i.e., to deal with the case
$b_i=\cdots=b_{i+\ell}$, we need another family of minimax methods, which involve the notion of genus.

\medskip

\BSS{\label[Subsection:2.b]. Symmetric mountain pass methods using genus}
In this section, we use an idea from Rabinowitz [\cite[R]] to define another family of minimax methods.
Here the notion of genus plays a role.

\medskip

\claim Definition.
Let $E$ be a Banach space.  For a closed set $A\subset E\setminus\{ 0\}$, which is symmetric with respect
to $0$, i.e., $-A=A$, we define $\g(A) = n$ if and only if there exists an odd map 
$\varphi\in C(A,\R^n\setminus \{ 0\})$ and $n$ is the smallest integer with this property.  When there 
are no odd map $\varphi\in C(A,\R^n\setminus \{ 0\})$ with this property for any $n\in\N$, we define $\g(A)=\infty$.  
Finally we set $\g(\emptyset)=0$.

\medskip

We refer to [\cite[R]] for fundamental properties of the genus.  

Our setting is different from [\cite[R]]; our functional is invariant under the following $\Z_2$-action:
    $$  \Z_2\times\RE \to \RE;\, (\pm 1,\lambda,u)\mapsto (\lambda,\pm u),      \eqno\label[2.6]
    $$
that is, $I(\lambda,-u)=I(\lambda,u)$.  
Remarking that there is no critical points in the $\Z_2$-invariants $\{ (\lambda,0);\, \lambda\in\R\}$,
we modify the arguments in [\cite[R]].

We define our second family of minimax sets as follows:
    $$  \eqalign{
        \Lambda_j &= \{ \gamma(\overline{D_{j+\ell}\setminus Y});\, \ell\geq 0,\, \gamma\in\Gamma_{j+\ell},\,
                \hbox{$Y\subset D_{j+\ell}\setminus\{ 0\}$ is closed,} \cr
            &\qquad \hbox{symmetric with respect to $0$ and $\g(Y)\leq \ell$} \},\cr
        c_j &= \inf_{A\in\Lambda_j} \max_{(\lambda,u)\in A} I(\lambda,u).\cr}
    $$
Here we summarize fundamental properties of $\Lambda_j$.  Here we use a projection\m $P_2:\, \RE\to\E$ 
defined by
    $$  P_2(\lambda,u)=u \quad \hbox{for}\ (\lambda,u)\in \RE.
    $$

\proclaim Lemma \label[Lemma:2.4].
\item{(i)} $\Lambda_j\not=\emptyset$ for all $j\in\N$.
\item{(ii)} $\Lambda_{j+1}\subset \Lambda_j$ for all $j\in\N$.
\item{(iii)} Let $\psi(\lambda,u)=(\psi_1(\lambda,u),\psi_2(\lambda,u)):\, \RE\to\RE$ be a continuous
map with properties
    $$  \eqalignno{
        &\psi_1(\lambda,-u)=\psi_1(\lambda,u), \ 
        \psi_2(\lambda,-u)=-\psi_2(\lambda,u) \ \hbox{for all}\ (\lambda,u)\in \RE, &\label[2.7]\cr
        &\psi(\lambda,u)=(\lambda,u) \quad \hbox{if}\ I(\lambda,u) \leq B_m-1.           &\label[2.8]\cr}
    $$
Then for $A\in \Lambda_j$, we have $\psi(A)\in \Lambda_j$.
\item{(iv)} For $A\in\Lambda_j$ and a closed set $Z$, which is invariant under $\Z_2$-action \ref[2.6],
i.e., $(\lambda,-u)\in Z$ for all $(\lambda,u)\in Z$, with $0\not\in \overline{P_2(Z)}$,
    $$  \overline{A\setminus Z} \in \Lambda_{j-i}, \ \hbox{where}\ i=\g(\overline{P_2(Z)}).
    $$
\item{(v)} $A\cap\partial \Omega_m\not=\emptyset$ for any $A\in \Lambda_j$.  Here $\Omega_m$ is defined
in \ref[2.3].

\medskip

\claim Proof.
(i), (ii) follow from the definition of $\Lambda_j$.   \m
(iii) Suppose $\psi(\lambda,u):\, \RE\to \RE$ satisfies \ref[2.7]--\ref[2.8].  Then it is easy to see
$\psi\circ \gamma\in \Gamma_j$ for all $\gamma\in\Gamma_j$.  Thus (iii) holds. \m
(iv) Following and modifying the argument in Sections \raw[7]--\raw[8] of [\cite[R]], we can show (iv).
For the sake of completeness, we give a proof in the Appendix. \m
(v) Suppose $A=\gamma(\overline{D_{j+\ell}\setminus Y})\in\Lambda_j$, where $\gamma\in\Gamma_{j+\ell}$,
$\g(Y)\leq\ell$, and let $U$ be the connected component of $\calO=\gamma^{-1}({\rm int}\, \Omega_m)$
containing $0$.  It is easy to see 
    $$  0\in U, \quad U\subset {\rm int}\, D_{j+\ell},
    $$
from which we have $\g(\partial U)=j+\ell$.  Thus
    $$  \g(\overline{\partial U\setminus Y})  \geq \g(\partial U)-\g(Y) \geq j.
    $$
In particular, $\overline{\partial U\setminus Y}\not=\emptyset$.
Since $\gamma(\overline{\partial U\setminus Y})
\subset A\cap\partial\Omega_m$, we have $A\cap\partial\Omega_m\not=\emptyset$.  \QED

\medskip

As fundamental properties of $c_n$, we have

\medskip

\proclaim Lemma \label[Lemma:2.5].
\item{(i)} $B_m\leq c_1\leq c_2\leq \cdots\leq c_j\leq c_{j+1}\leq\cdots$.
\item{(ii)} $c_j\leq b_j$ for all $j\in \N$.

\medskip

\claim Proof.
(i) By (v) of Lemma \ref[Lemma:2.4], we have for any $A\in\Lambda_j$,
    $$  \max_{(\lambda,u)\in A} I(\lambda,u) \geq \inf_{(\lambda,u)\in\partial\Omega_m} I(\lambda,u) = B_m,
    $$
which implies $c_j\geq B_m$ for all $j\in\N$.  (ii) of Lemma \ref[Lemma:2.4] implies $c_j \leq c_{j+1}$.  \m
(ii) It is easy to see $\gamma(D_j)\in \Lambda_j$ for any $\gamma\in\Gamma_j$.  Thus we have
$c_j\leq b_j$.  \QED

\medskip

In the following section, we use a special deformation lemma to show $c_j$ ($j=\m 1,2,\cdots,k$) are 
attained by critical points.

\medskip

\BS{\label[Section:3]. Deformation argument and existence of critical points}
In this section we introduce a deformation result for $I(\lambda,u)$ and we show that
$c_j$ ($j=1,2,\cdots,k$) introduced in the previous section are achieved.

\BSScap{\label[Subsection:3.a].  Deformation result for $I(\lambda,u)$}{\raw[3.a].  Deformation result for I}
For $b\in \R$ we set
    $$  K_b = \{(\lambda,u)\in\RE;\, I(\lambda,u)=b,\, \partial_\lambda I(\lambda,u)=0,\,
        \partial_u I(\lambda,u)=0,\, P(\lambda,u)=0 \}.      \eqno\label[3.1]
    $$
Here $P(\lambda,u)$ is introduced in \ref[1.5].  
We note that $\partial_u I(\lambda,u)=0$ implies $P(\lambda,u)=0$.
We also use the following notation:
    $$  [I\leq c] =\{(\lambda,u)\in\E;\, I(\lambda,u)\leq c\} \quad \hbox{for}\ c\in\R.
    $$
We have the following deformation result.

\medskip

\proclaim Proposition \label[Proposition:3.1]. 
Assume \cond[g:1]--\cond[g:5] and $b<0$.  Then
\item{(i)} $K_b$ is compact in $\RE$ and $K_b\cap(\R\times\{0\})=\emptyset$.
\item{(ii)} For any open neighborhood $\calO$ of $K_b$ and $\overline\epsilon>0$ there exist
$\epsilon\in (0,\overline\epsilon)$ and a continuous map $\eta(t,\lambda,u):\, [0,1]\times\RE
\to\RE$ such that
\itemitem{\maru 1} $\eta(0,\lambda,u)=(\lambda,u)$ for all $(\lambda,u)\in\RE$.
\itemitem{\maru 2} $\eta(t,\lambda,u)=(\lambda,u)$ if  $(\lambda,u)\in [I\leq b-\overline\epsilon]$.
\itemitem{\maru 3} $I(\eta(t,\lambda,u))\leq I(\lambda,u)$ for all $(t,\lambda,u)\in [0,1]\times\RE$.
\itemitem{\maru 4} $\eta(1,[I\leq b+\epsilon]\setminus\calO)\subset [I\leq b-\epsilon]$, \m
        $\eta(1,[I\leq b+\epsilon])\subset [I\leq b-\epsilon]\cup\calO$.
\itemitem{\maru 5} If $K_b=\emptyset$, $\eta(1,[I\leq b+\epsilon])\subset [I\leq b-\epsilon]$. 
\itemitem{\maru 6} Writing $\eta(t,\lambda,u)=(\eta_1(t,\lambda,u),\eta_2(t,\lambda,u))$, we have
    $$  \eta_1(t,\lambda,-u)=\eta_1(t,\lambda,u), \quad \eta_2(t,\lambda,-u)=-\eta_2(t,\lambda,u)
    $$
for all $(t,\lambda,u)\in [0,1]\times \RE$.
\medskip

\noindent
Such a deformation result is usually obtained under the Palais-Smale compactness condition.
However it seems difficult to verify the standard Palais-Smale condition under \cond[g:1]--\cond[g:4].
In Section \ref[Section:4], we introduce a new version $(PSP)$ of Palais-Smale condition and we develop a new
deformation argument to prove Proposition \ref[Proposition:3.1].  We postpone a proof of Proposition \ref[Proposition:3.1] until
Section \ref[Section:4] and in this section we show $c_j$ ($j=1,2,\cdots,k$) are attained.

\medskip

\BSS{\label[Subsection:3.b]. Existence of critical points}
As an application of our Proposition \ref[Proposition:3.1] we show the following

\medskip

\claim Proposition \label[Proposition:3.2].
\item{(i)} For $j=1,2,\cdots,k$, $c_j<0$ and $c_j$ is a critical value of $I(\lambda,u)$.
\item{(ii)} If $c_j=c_{j+1}=\cdots=c_{j+q}=b<0$ ($j+q\leq k$), then
    $$  \g(P_2(K_b)) \geq q+1.
    $$
In particular, $\#(K_b)=\infty$ if $q\geq 1$.

\claim Proof.
$c_j<0$ ($j=1,2,\cdots,k$) follow from Proposition \ref[Proposition:2.2] and 
Lemma \ref[Lemma:2.5] (ii).  The argument for the fact that $K_{c_j}\not=\emptyset$ is similar to 
the proof of (ii).  So we omit it. \m
(ii) 
Suppose that $c_j=c_{j+1}=\cdots=c_{j+q}=b<0$.  Since $K_b$ is compact and 
$K_b\cap(\R\times\{0\})=\emptyset$, the projection $P_2(K_b)$ of $K_b$ onto $\E$ is compact,
symmetric with respect to $0$ and $0\not\in P_2(K_b)$.  Thus by the fundamental property 
of the genus,
    $$  \eqalign{
        &\g(P_2(K_b)) <\infty,\cr
        &\hbox{there exists $\delta>0$ small such that $\g(\overline{P_2(N_\delta(K_b))})=\g(P_2(K_b))$}. \cr}
    $$
Here we denote $\delta$-neighborhood of a set $A\subset\RE$ by $N_\delta(A)$, i.e.,
    $$  N_\delta(A) =\{ (\lambda,u);\, \dist((\lambda,u), A)\leq \delta\},
    $$
where 
    $$  \dist((\lambda,u), A)
        =\inf_{(\lambda',u')\in A} \sqrt{\abs{\lambda-\lambda'}^2+\norm{u-u'}_{H^1}^2}.
    $$
By Proposition \ref[Proposition:3.1], there exist $\epsilon>0$ small and $\eta:\, [0,1]\times\RE\to\RE$ such that
    $$  \eqalign{
        &\eta(1,[I\leq b+\epsilon]\setminus N_\delta(K_b)) \subset [I\leq b-\epsilon],\cr
        &\eta(t,\lambda,u)=(\lambda,u) \quad \hbox{if}\ I(\lambda,u)\leq b-\half. \cr}
    $$
We note that $B_m-1\leq b-\half$.\m
We take $A\in \Lambda_{j+q}$ such that $A\subset [I\leq b+\epsilon]$.  Then
    $$  \eta(1,\overline{A\setminus N_\delta(K_b)}) \subset [I\leq b-\epsilon].  \eqno\label[3.2]
    $$
If $\g(P_2(K_b))\leq q$, we have $\g(\overline{P_2(N_\delta(K_b))})\leq q$.  By (iv) of Lemma \ref[Lemma:2.4],
    $$  \overline{A\setminus N_\delta(K_b)} \in \Lambda_j.      \eqno\label[3.3]
    $$
\ref[3.2] and \ref[3.3] imply $c_j\leq b-\epsilon$, which is a contradiction.  Thus 
$\g(P_2(K_b))\geq q+1$.  \QED

\medskip

Now we can show

\medskip

\claim Proof of (i) of Theorem \ref[Theorem:0.2].
Clearly (i) of Theorem \ref[Theorem:0.2] follows from Proposition \ref[Proposition:3.2].  \QED

\medskip

\claim Proof of (ii) of Theorem \ref[Theorem:0.2].
Under the condition \ref[0.2], we have $m_k=0$ for all $k\in\N$.  Thus we have $c_j\leq b_j<0$ for all
$j\in\N$ and $c_j$ ($j\in\N$) are critical values of $I(\lambda,u)$.  We need to show $c_j\to 0$
as $j\to\infty$.

Arguing indirectly, we assume $c_j\to \overline c<0$ as $j\to\infty$.  
Then $K_{\overline c}$ is compact and $K_{\overline c}\cap(\R\times\{0\})=\emptyset$.  Set
    $$  q=\g(P_2(K_{\overline c})) <\infty
    $$
and choose $\delta>0$ small such that
    $$  \g(\overline{P_2(N_\delta(K_{\overline c}))}) = \g(P_2(K_{\overline c})) = q.
    $$
As in the proof of Proposition \ref[Proposition:3.2], there exist $\epsilon>0$ small and $\eta:\,[0,1]\times\RE\to \RE$
such that
    $$  \eqalignno{
        &\eta(1,[I\leq {\overline c}+\epsilon]\setminus N_\delta(K_{\overline c})) 
                                    \subset [I\leq {\overline c}-\epsilon], &\label[3.4]\cr
        &\eta(t,\lambda,u)=(\lambda,u) \quad \hbox{if}\ I(\lambda,u) \leq B_m-1. \cr}
    $$
We choose $j\gg 1$ so that $c_j \geq {\overline c}-\epsilon$ and take $A\in \Lambda_j$ such that
$A\subset [I\leq {\overline c}+\epsilon]$.  Then we have
    $$  \overline{A\setminus N_\delta(K_{\overline c})} \in \Lambda_{j-q}.      \eqno\label[3.5]
    $$
\ref[3.4] and \ref[3.5] imply $c_{j-q}\leq {\overline c}-\epsilon$.  Since we can take $j$ arbitrary large, we have
$\lim_{j\to\infty} c_j \leq {\overline c}-\epsilon$.  This is a contradiction.  \QED

\medskip

\BScap{\label[Section:4]. $(PSP)$ condition and construction of a flow}{4. PSP condition and construction of a flow}
In this section we give a new type of deformation argument for our functional $I(\lambda,u)$. 
Our deformation argument is inspired by our previous work [\cite[HIT]].

\BSScap{\label[Subsection:4.a]. $(PSP)$ condition}{\raw[4.a]. PSP condition}
Since it is difficult to verify the standard Palais-Smale condition for $I(\lambda,u)$ under 
the conditions \cond[g:1]--\cond[g:5], we introduce a new type of Palais-Smale condition
$(PSP)_b$, which is weaker than the standard Palais-Smale condition and which takes the scaling
property of $I(\lambda,u)$ into consideration through the Pohozaev functional $P(\lambda,u)$ .

\medskip

\claim Definition.
For $b\in\R$, we say that $I(\lambda,u)$ satisfies $(PSP)_b$ condition if and only if the following
holds.
\itemitem{$(PSP)_b$} If a sequence $(\lambda_n,u_n)_{n=1}^\infty\subset \RE$ satisfies as
$n\to\infty$
    $$  \eqalignno{
        &I(\lambda_n,u_n) \to b,                            &\label[4.1]\cr
        &\partial_\lambda I(\lambda_n,u_n) \to 0,           &\label[4.2]\cr
        &\partial_u I(\lambda_n,u_n) \to 0 \quad \hbox{strongly in}\ (\E)^*,            &\label[4.3]\cr
        &P(\lambda_n,u_n) \to 0,                            &\label[4.4]\cr}
    $$
then $(\lambda_n,u_n)_{n=1}^\infty$ has a strongly convergent subsequence in $\RE$.

\medskip

\noindent
First we observe that $(PSP)_b$ holds for $I(\lambda,u)$ for $b<0$.

\medskip

\proclaim Proposition \label[Proposition:4.1].
Assume \cond[g:1]--\cond[g:4].  Then $I(\lambda,u)$ satisfies $(PSP)_b$ for $b<0$.

\medskip

\claim Proof.
Let $b<0$ and suppose that $(\lambda_n,u_n)_{n=1}^\infty$ satisfies \ref[4.1]--\ref[4.4]. 
We will show that $(\lambda_n,u_n)_{n=1}^\infty$ has a strongly convergent subsequence.
Proof consists of several steps.

\smallskip

\noindent
{\sl Step 1: $\lambda_n$ is bounded from below as $n\to\infty$.}

\smallskip

\noindent
Since
    $$  P(\lambda_n,u_n) =N(I(\lambda_n,u_n)+{m\over 2}e^{\lambda_n}) - \norm{\nabla u_n}_2^2,
    $$
we have from \ref[4.1], \ref[4.4] that
    $$  {m\over 2}\liminf_{n\to \infty} e^{\lambda_n} \geq - b>0.
    $$
Thus $\lambda_n$ is bounded from below as $n\to \infty$.

\smallskip

\noindent
{\sl Step 2: $\norm{u_n}_2^2 \to m$ as $n\to \infty$.}

\smallskip

\noindent
Since $\partial_\lambda I(\lambda_n,u_n)={e^{\lambda_n}\over 2}(\norm{u_n}_2^2-m)$, 
it follows from \ref[4.2] and Step 1 that $\norm{u_n}_2^2\to m$.

\smallskip

\noindent
{\sl Step 3: $\norm{\nabla u_n}_2^2$ and $\lambda_n$ are bounded as $n\to\infty$.}

\smallskip

\noindent
We have
    $$  \partial_u I(\lambda_n,u_n)u_n = \norm{\nabla u_n}_2^2 -\intRN g(u_n)u_n + e^{\lambda_n}\norm{u_n}_2^2.
                                                                \eqno\label[4.5]
    $$
By \cond[g:2] and \cond[g:3], for any $\delta>0$ there exists $C_\delta>0$ such that
    $$  \abs{g(\xi)\xi} \leq C_\delta\abs\xi^2 + \delta\abs\xi^{p+1}    \quad \hbox{for all}\ \xi\in\R.
    $$
Thus
    $$  \pabs{\intRN g(u)u} \leq C_\delta\norm u_2^2 + \delta\norm u_{p+1}^{p+1}
        \quad \hbox{for all}\ u\in \E.
    $$
Since $p=1+{4\over N}$, by Gagliardo-Nirenberg inequality there exists $C_N>0$ such that
    $$  \norm u_{p+1}^{p+1} \leq C_N\norm{\nabla u}_2^2\norm u_2^{p-1}
        \quad \hbox{for all}\ u\in \E.
    $$
Thus it follows from \ref[4.5] that
    $$  \norm{\nabla u_n}_2^2 -C_\delta\norm{u_n}_2^2 -\delta C_N\norm{\nabla u_n}_2^2\norm{u_n}_2^{p-1}
        + e^{\lambda_n}\norm{u_n}_2^2
        \leq \epsilon_n\sqrt{\norm{\nabla u_n}_2^2 + \norm{u_n}_2^2},
    $$
where $\epsilon_n = \norm{\partial_u I(\lambda_n,u_n)}_{(\E)^*}\to 0$.

By Step 2,
    $$  \left(1-\delta C_N(m+o(1))^{p-1\over 2}\right)\norm{\nabla u_n}_2^2 
        +(e^{\lambda_n}- C_\delta)(m+o(1)) 
        \leq \epsilon_n \sqrt{\norm{\nabla u_n}_2^2 + m + o(1)}.
    $$
Choosing $\delta >0$ small so that $\delta C_Nm^{p-1\over 2}<\half$, we observe that $\norm{\nabla u_n}_2^2$
and $e^{\lambda_n}$ are bounded as $n\to\infty$.

\smallskip

\noindent
{\sl Step 4: Conclusion}

\smallskip

\noindent
By Steps 1--3, $(\lambda_n,u_n)_{n=1}^\infty$ is a bounded sequence in $\RE$.  After extracting a subsequence
--- still denoted by $(\lambda_n,u_n)_{n=1}^\infty$ ---,  we may assume that $\lambda_n\to\lambda_0$ and
$u_n\wlimit u_0$ weakly in $\E$ for some $(\lambda_0,u_0)\in\RE$.  
By \cond[g:2], \cond[g:3], we have
    $$ \intRN g(u_n)u_0 \to \intRN g(u_0)u_0, \quad \intRN g(u_n)u_n \to \intRN g(u_0)u_0.
    $$
Thus, we deduce from $\partial_u I(\lambda_n,u_n)u_n\to 0$ and 
$\partial_u I(\lambda_n,u_n)u_0\to 0$ that
    $$  \norm{\nabla u_n}_2^2 +e^{\lambda_0}\norm{u_n}_2^2
        \to \norm{\nabla u_0}_2^2 +e^{\lambda_0}\norm{u_0}_2^2,
    $$
which implies $u_n\to u_0$ strongly in $\E$.                    \QED

\medskip

\claim Remark \label[Remark:4.2].  
For $b=0$, $(PSP)_0$ does not hold for $I(\lambda,u)$.  In fact, for a sequence $(\lambda_n,0)_{n=1}^\infty$ with
$\lambda_n\to-\infty$, we have
    $$  \eqalign{
        &I(\lambda_n,0)=-{e^{\lambda_n}\over 2}m \to 0,\quad \partial_\lambda I(\lambda_n,0) =-{e^{\lambda_n}\over 2}m \to 0,\cr
        &\partial_u I(\lambda_n,0)=0, \qquad\qquad\quad   P(\lambda_n,0)=0. \cr}
    $$
But  $(\lambda_n,0)_{n=1}^\infty$ has no convergent subsequences.

\medskip

As a corollary to Proposition \ref[Proposition:4.1], we have

\medskip

\proclaim Corollary \label[Corollary:4.3].
For $b<0$, $K_b$ defined in \ref[3.1] is compact in $\RE$ and satisfies $K_b\cap(\R\times\{0\})=\emptyset$. 

\medskip

\claim Proof.
$K_b$ is compact since $I(\lambda,u)$ satisfies $(PSP)_b$.  $K_b\cap (\R\times\{ 0\})\not=\emptyset$
follows from the fact that $\partial_\lambda I(\lambda,0)=-{e^\lambda\over 2}m\not= 0$.
\QED

\medskip


\BSScap{\label[Subsection:4.b]. Functional $J(\tlu)$}{\raw[4.b]. Functional J}
To construct a deformation flow we need an augmented functional\m 
$J(\theta,\lambda,u):\,\RRE\to \R$ defined by
    $$  J(\tlu)= \half e^{(N-2)\theta}\norm{\nabla u}_2^2 -e^{N\theta}\intRN G(u)
        +{e^\lambda\over 2}\left( e^{N\theta}\norm u_2^2 -m\right).
    $$
We introduce $J(\tlu)$ to make use of the scaling property of $I(\lambda,u)$.  
As a basic property of $J(\theta,\lambda,u)$ we have
    $$  I(\lambda,u(x/e^\theta)) = J(\tlu) \quad \hbox{for all}\ (\tlu)\in \RRE.
                                                \eqno\label[4.6]
    $$
We will construct our deformation flow for $I(\lambda,u)$ through a deformation flow for $J(\tlu)$.

$J(\tlu)$ satisfies the following properties.

\medskip

\proclaim Lemma \label[Lemma:4.4].
For all $(\tlu)\in \RRE$, $h\in\E$ and $\beta\in\R$,
    $$  \eqalignno{
        &\partial_\theta J(\tlu(x)) = P(\lambda,u(x/e^\theta)),                 &\label[4.7]\cr
        &\partial_\lambda J(\tlu(x)) = \partial_\lambda I(\lambda,u(x/e^\theta)),   &\label[4.8]\cr
        &\partial_u J(\tlu(x))h(x) = \partial_u I(\lambda,u(x/e^\theta)) h(x/e^\theta), 
                                                                                &\label[4.9]\cr
        &J(\theta+\beta,\lambda,u(e^\beta x)) = J(\tlu(x)).                     &\label[4.10]\cr}
    $$

\claim Proof.
We compute that
    $$  \eqalign{
    \partial_\theta J(\tlu(x)) &= {N-2\over 2}e^{(N-2)\theta}\norm{\nabla u}_2^2
        + Ne^{N\theta}\left({e^\lambda\over 2}\norm u_2^2 -\intRN G(u)\right) \cr
    &= {N-2\over 2} \norm{\nabla(u(x/e^\theta x))}_2^2
        + N \left({e^\lambda\over 2}\norm{u(x/e^\theta)}_2^2 -\intRN G(u(x/e^\theta))\right) \cr
    &= P(\lambda,u(x/e^\theta)),\cr
    \partial_\lambda J(\tlu(x)) &= {e^\lambda\over 2}\left(e^{N\theta}\norm u_2^2 -m\right)
     = {e^\lambda\over 2}\left(\norm{u(x/e^\theta)}_2^2 -m\right)\cr
    &= \partial_\lambda I(\lambda,u(x/e^\theta)),\cr
    \partial_u J(\tlu(x)) h(x) &= e^{(N-2)\theta}(\nabla u,\nabla h)_2 +e^\lambda e^{N\theta}(u,h)_2
        - e^{N\theta}\intRN g(u(x))h(x) \cr
    &= (\nabla u(x/e^\theta),\nabla h(x/e^\theta))_2 +e^\lambda (u(x/e^\theta),h(x/e^\theta))_2\cr
    &\qquad - \intRN g(u(x/e^\theta))h(x/e^\theta) \cr
    &= \partial_u I(\lambda,u(x/e^\theta)) h(x/e^\theta). \cr}
    $$
Thus we have \ref[4.7]--\ref[4.9].  \ref[4.10] follows from \ref[4.6].  \QED

\medskip

To analyze $J(\tlu)$, it is natural to regard $\RRE$ as a Hilbert manifold with a metric
related to \ref[4.6].  More precisely, we write $M=\RRE$.  We note that
    $$  T_{(\tlu)}M= \RRE \quad \hbox{for}\ (\tlu)\in M
    $$
and we introduce a metric $\langle\cdot,\cdot\rangle_{(\tlu)}$ on $T_{(\tlu)} M$ by
    $$  \eqalign{
        &\langle(\alpha,\nu,h), (\alpha',\nu',h')\rangle_{(\tlu)}
        = \alpha\alpha' + \nu\nu' + e^{(N-2)\theta} (\nabla h,\nabla h')_2
            + e^{N\theta} (h,h')_2,\cr
        &\norm{(\alpha,\nu,h(x))}_{(\tlu)} 
        =\sqrt{\langle(\alpha,\nu,h),(\alpha,\nu,h)\rangle_{(\tlu)}} \cr}
    $$
for $(\alpha,\nu,h)$, $(\alpha',\nu',h')\in T_{(\tlu)}M$.
We also denote the dual norm of $\norm\cdot_{(\tlu)}$ by $\norm\cdot_{(\tlu),*}$, that is,
    $$  \norm f_{(\tlu),*} = \sup_{\norm{(\alpha,\nu,h)}_{(\tlu)} \leq 1} \abs{f(\alpha,\nu,h)}
        \quad \hbox{for}\ f\in T_{(\tlu)}^*(M).
                                                    \eqno\label[4.11]
    $$
It is easily seen that $(M,\langle\cdot,\cdot\rangle)$ is a complete Hilbert manifold.  We note
that $\langle\cdot,\cdot\rangle_{(\tlu)}$ and $\norm\cdot_{(\tlu)}$ depend only on $\theta$. 
So sometimes we denote them by $\langle\cdot,\cdot\rangle_{(\theta,\cdot,\cdot)}$, 
$\norm\cdot_{(\theta,\cdot,\cdot)}$.  We have
    $$  \eqalignno{
    \norm{(\alpha,\nu,h)}_{(\theta,\cdot,\cdot)} &= \alpha^2+\nu^2 
        +e^{(N-2)\theta}\norm{\nabla h}_2^2 +e^{N\theta}\norm{h}_2^2 \cr
        &= \alpha^2+\nu^2 + \norm{h(x/e^\theta)}_{H^1}^2 \cr
        &= \norm{(\alpha,\nu,h(x/e^\theta))}_{(0,\cdot,\cdot)}
            \quad \hbox{for}\ (\alpha,\nu,h)\in T_{(\tlu)}M.     &\label[4.12]\cr}
    $$
We also have for all $(\alpha,\nu,h)\in T_{(\theta,\cdot,\cdot)}M$ and $\beta\in\R$
    $$  \norm{(\alpha,\nu,h(e^\beta x))}_{(\theta+\beta,\cdot,\cdot)}
        = \norm{(\alpha,\nu,h(x))}_{(\theta,\cdot,\cdot)}.          \eqno\label[4.13]
    $$
We also define a distance $\distM(\cdot,\cdot)$ on $M$ by
    $$  \eqalignno{
        &\distM((\theta_0,\lambda_0,h_0),(\theta_1,\lambda_1,h_1)) \cr
        =& \inf\Bigl\{\int_0^1 \norm{\dot\sigma(t)}_{\sigma(t)} \, dt;\,
            \sigma(t)\in C^1([0,1],M),\, \cr
        &\qquad \sigma(0)=(\theta_0,\lambda_0,h_0),\, 
            \sigma(1)=(\theta_1,\lambda_1,h_1)\Bigr\}.          &\label[4.14]\cr}
    $$
By the property \ref[4.13], we have for all $\beta\in\R$
    $$  \distM((\theta_0+\beta,\lambda_0,u_0(e^\beta x)),(\theta_1+\beta,\lambda_1,u_1(e^\beta x)))
        = \distM((\theta_0,\lambda_0,u_0(x)),(\theta_1,\lambda_1,u_1(x))).
                                                                \eqno\label[4.15]
    $$
Using notation
    $$  \Diff =(\partial_\theta,\partial_\lambda,\partial_u),
    $$
we have

\medskip

\proclaim Lemma \label[Lemma:4.5].  
For $(\tlu)\in M$, we have
    $$  \eqalign{
        &\norm{\Diff J(\tlu)}_{(\tlu),*}   \cr
        =&\left( \abs{P(\lambda,u(x/e^\theta))}^2
            + \abs{\partial_\lambda I(\lambda,u(x/e^\theta))}^2 
            + \norm{\partial_u I(\lambda,u(x/e^\theta))}_{(\E)^*}^2 \right)^{1/2}.\cr}
    $$

\claim Proof.
By Lemma \ref[Lemma:4.4], we have
    $$  \eqalign{
        &\Diff J(\tlu)(\alpha,\nu,h) \cr
        =& P(\lambda,u(x/e^\theta))\alpha + \partial_\lambda I(\lambda,u(x/e^\theta))\nu
            + \partial_u I(\lambda,u(x/e^\theta)) h(x/e^\theta). \cr}
    $$
Noting \ref[4.12], the conclusion of Lemma \ref[Lemma:4.5] follows from the definition \ref[4.11].  \QED

\medskip

For $b\in \R$, we use notation
    $$  \wKb =\{ (\tlu)\in M;\, J(\tlu)=b, \, \Diff J(\tlu)=(0,0,0)\}.
    $$
By \ref[4.6]--\ref[4.9], we observe that
    $$  \wKb =\{(\theta,\lambda,u(e^\theta x));\, \theta\in\R,\, (\lambda,u)\in K_b\}.
    $$
We also use notation for $(\tlu)\in M$ and $\widetilde A\subset M$
    $$  \distM((\tlu),\widetilde A) 
        = \inf_{(\theta',\lambda',u')\in \widetilde A}\distM((\tlu),(\theta',\lambda',u')).
    $$
From $(PSP)_b$ condition for $I(\lambda,u)$, we deduce the following

\medskip

\proclaim Proposition \label[Proposition:4.6].
For $b<0$, $J(\tlu)$ satisfying the following property:
\itemitem{$(\widetilde{PSP})_b$}
For any sequence $(\theta_n,\lambda_n,u_n)_{n=1}^\infty\subset M$ with
    $$  \eqalignno{
        &J(\theta_n,\lambda_n,u_n) \to b,                                           &\label[4.16]\cr
        &\norm{\Diff J(\theta_n,\lambda_n,u_n)}_{(\theta_n,\lambda_n,u_n),*}\to 0
                \quad \hbox{as}\ n\to\infty,    &\label[4.17]\cr}
    $$
we have
    $$  \distM((\theta_n,\lambda_n,u_n), \wKb) \to 0.               \eqno\label[4.18]
    $$

\claim Proof.
Suppose that $(\theta_n,\lambda_n,u_n)_{n=1}^\infty$ satisfies \ref[4.16]--\ref[4.17].
It suffices to show that $(\theta_n,\lambda_n,u_n)_{n=1}^\infty$ has a subsequence with the property
\ref[4.18].

Setting $\hat u_n(x)=u_n(x/e^{\theta_n})$, we have by Lemma \ref[Lemma:4.5] that
    $$  \eqalign{
        &I(\lambda_n,\hat u_n)\to b <0, \cr
        &P(\lambda_n,\hat u_n)\to 0, \quad \partial_\lambda I(\lambda_n,\hat u_n)\to 0, \quad
        \partial_u I(\lambda_n,\hat u_n)\to 0 \ \hbox{strongly in}\ (\E)^*. \cr}
    $$
Thus by Proposition \ref[Proposition:4.1], there exists a subsequence --- still denoted by $(\lambda_n,\hat u_n)_{n=1}^\infty$
--- and $(\lambda_0,\hat u_0)\in \RE$ such that
    $$  \lambda_n \to \lambda_0 \quad \hbox{and}\quad \hat u_n\to \hat u_0 \ \hbox{strongly in}\ \E.
    $$
Note that $(\lambda_0,\hat u_0)\in K_b$ and thus $(\theta_n,\lambda_0, \hat u_0(e^{\theta_n} x))
\in \wKb$.  By \ref[4.15], we have
    $$  \eqalign{
        &\distM((\theta_n,\lambda_n,u_n), \wKb) 
        \leq \distM((\theta_n,\lambda_n,u_n), (\theta_n,\lambda_0, \hat u_0(e^{\theta_n} x))) \cr
        = & \distM((0,\lambda_n,\hat u_n), (0,\lambda_0,\hat u_0(x))
        \leq \left(\abs{\lambda_n-\lambda_0}^2 + \norm{\hat u_n-\hat u_0}_{H^1}^2 \right)^{1/2}
        \to 0 \quad \hbox{as}\ n\to\infty.  \cr}
    $$
\QED

\medskip

As a corollary to Proposition \ref[Proposition:4.6], we have the following uniform estimate of $\Diff J(\tlu)$ outside
of $\rho$-neighborhood of $\wKb$.

\medskip

\proclaim Corollary \label[Corollary:4.7].
Assume $b<0$.  Then for any $\rho>0$ there exists $\delta_\rho>0$ such that for $(\tlu)\in M$
    $$  \abs{J(\tlu)-b} < \delta_\rho \quad \hbox{and}\quad
        \distM((\tlu),\wKb) \geq \rho
    $$
imply
    $$  \norm{\Diff J(\tlu)}_{(\tlu),*} \geq \delta_\rho.                   \QED
    $$

We remark that $\wKb$ is not compact in $M$ but Corollary \ref[Corollary:4.7] gives us a uniform lower bound of
$\norm{\Diff J(\tlu)}_{(\tlu),*}$ outside of $\rho$-neighborhood of $\wKb$, which enables us
to construct a deformation flow for $J(\tlu)$.

\medskip


\BSScap{\label[Subsection:4.c]. Deformation flow for $J(\tlu)$}{\raw[4.c]. Deformation flow for J}
In this section we give a deformation result for $J(\tlu)$.  We need the following
notation:
    $$  \eqalign{
    &[J\leq c]_M =\{ (\tlu)\in M;\, J(\tlu)\leq c\} \quad \hbox{for}\ c\in\R,\cr
    &\wN_\rho(\widetilde A) = \{ (\tlu)\in M;\, \distM((\tlu), \widetilde A) \leq \rho\} \quad 
        \hbox{for}\ \widetilde A\subset M \ \hbox{and}\ \rho>0. \cr}
    $$
We have the following deformation result.

\medskip

\proclaim Proposition \label[Proposition:4.8].
Assume $b<0$.  Then for any $\overline\epsilon>0$ and $\rho>0$ there exist 
$\epsilon\in (0,\overline\epsilon)$ and a continuous map $\weta(t,\tlu):\, [0,1]\times M\to M$ 
such that
\item{\maru 1} $\weta(0,\tlu)=(\tlu)$ for all $(\tlu)\in M$.
\item{\maru 2} $\weta(t,\tlu)=(\tlu)$ if  $(\tlu)\in [J\leq b-\overline\epsilon]_M$.
\item{\maru 3} $J(\weta(t,\tlu))\leq J(\tlu)$ for all $(t,\tlu)\in [0,1]\times M$.
\item{\maru 4} $\weta(1,[J\leq b+\epsilon]_M\setminus \wN_\rho(\wKb))\subset [J\leq b-\epsilon]_M$, \m
        $\weta(1,[J\leq b+\epsilon]_M)\subset [J\leq b-\epsilon]_M\cup \wN_\rho(\wKb)$.
\item{\maru 5} If $K_b=\emptyset$, $\weta(1,[J\leq b+\epsilon]_M)\subset [J\leq b-\epsilon]_M$. 
\item{\maru 6} We write $\weta(t,\tlu)=(\weta_0(t,\tlu), \weta_1(t,\tlu),\weta_2(t,\tlu))$.
Then $\weta_0(t,\tlu)$, $\weta_1(t,\tlu)$ are even in $u$ and $\weta_2(t,\tlu)$ is odd in $u$.  
That is, for all  $(t,\tlu)\in [0,1]\times M$
    $$  \eqalign{
        &\weta_0(t,\theta, \lambda,-u)=\weta_0(t,\tlu), \quad 
        \weta_1(t,\theta,\lambda,-u)=\weta_1(t,\tlu), \cr
        &\weta_2(t,\theta,\lambda,-u)=-\weta_2(t,\tlu).  \cr}
    $$

\claim Proof.
Let $M'=\{ (\tlu)\in M;\, \Diff J(\tlu) \not= (0,0,0)\}$.  It is well-known that there
exists a pseudo-gradient vector field $\calV:\, M'\to TM$ such that for 
$(\tlu)\in M'$

\smallskip

\itemitem{(1)} $\norm{\calV(\tlu)}_{(\tlu)} \leq 2\norm{\Diff J(\tlu)}_{(\tlu),*}$,
\itemitem{(2)} $\Diff J(\tlu)\calV(\tlu) \geq \norm{\Diff J(\tlu)}_{(\tlu),*}^2$,
\itemitem{(3)} $\calV:\, M'\to \RRE$ is locally Lipschitz continuous.

\smallskip

\noindent 
We can also have
\itemitem{(4)} $\calV(\tlu)=(\calV_0(\tlu),\calV_1(\tlu),\calV_2(\tlu))$ satisfies
    $$  \eqalign{
        &\quad  \calV_0(\theta,\lambda,-u) = \calV_0(\tlu), \quad 
                \calV_1(\theta,\lambda,-u) = \calV_1(\tlu), \cr
        &\quad  \calV_2(\theta,\lambda,-u) = -\calV_2(\tlu). \cr}
    $$

\noindent
For a given $\rho>0$ we choose $\delta_\rho>0$ by Corollary \ref[Corollary:4.7] so that
    $$  \abs{J(\tlu)-b} < \delta_\rho \quad \hbox{and}\quad
        (\tlu) \not\in \wN_{\rho/3}(\wKb) \quad \hbox{imply} \quad
        \norm{\Diff J(\tlu)}_{(\tlu),*}\geq \delta_\rho.        \eqno\label[4.19]
    $$
We choose a locally Lipschitz continuous function $\varphi:\, M\to [0,1]$ such that
    $$  \eqalign{
        &\varphi(\tlu) = 1 \quad \hbox{for} \ (\tlu)\in M\setminus\wN_{{2\over 3}\rho}(\wKb),
                                                                \cr
        &\varphi(\tlu) = 0 \quad \hbox{for} \ (\tlu)\in \wN_{{1\over 3}\rho}(\wKb), \cr
        &\varphi(\theta,\lambda,-u)=\varphi(\tlu) \quad \hbox{for all}\ (\tlu)\in M. \cr}
    $$
For $\overline\epsilon>0$ we may assume $\overline\epsilon\in (0,\delta_\rho)$ and
we choose a locally Lipschitz continuous function $\psi:\, \R\to [0,1]$ such that
    $$  \psi(s) = \cases{
        1   &for $s\in [b-{\overline\epsilon\over 2}, b+{\overline\epsilon\over 2}]$, \cr
        0   &for $s\in \R\setminus [b-\overline\epsilon, b+\overline\epsilon]$. \cr}
    $$
We consider the following ODE in $M$:
    $$  \left\{\ \eqalign{
        &{d\weta\over dt} 
        = -\varphi(\weta)\psi(J(\weta)) {\calV(\weta)\over \norm{\calV(\weta)}_{\weta}}, \cr
        &\weta(0,\tlu)=(\tlu).  \cr}\right.
    $$
For $\epsilon\in (0,\overline\epsilon)$ small, $\weta(t,\tlu)$ has the desired properties
\maru 1--\maru 6.  We show just the first part of \maru 4:
    $$  \weta(1,[J\leq b+\epsilon]_M\setminus \wN_\rho(\wKb))\subset [J\leq b-\epsilon]_M.
                                                                \eqno\label[4.20]
    $$
We can check properties \maru 1--\maru 3 easily and we use them in what follows.  We also
note that 
    $$  \norm{{d\weta\over dt}(t)}_{\weta(t)} \leq 1 \quad \hbox{for all}\  t.
                                                                \eqno\label[4.21]
    $$
For $\epsilon\in (0,{\overline\epsilon\over 2})$, which we choose later, we assume $\weta(t)=
\weta(t,\tlu)$ satisfies
    $$  \weta(0)\in [J\leq b+\epsilon]_M \setminus \wN_\rho(\wKb).
    $$
If $\weta(1)\not\in [J\leq b-\epsilon]_M$, we have $J(\weta(t))\in [b-\epsilon,b+\epsilon]$ for all
$t\in [0,1]$.  We consider 2 cases:

\smallskip

\itemitem{Case 1:} $\weta(t) \not\in \wN_{{2\over 3}\rho}(\wKb)$ for all $t\in [0,1]$,
\itemitem{Case 2:} $\weta(t_0) \in \wN_{{2\over 3}\rho}(\wKb)$ for some $t_0\in [0,1]$.

\smallskip

\noindent
First we consider Case 1.  By \ref[4.19] we have 
    $$  \norm{\Diff J(\weta(t))}_{\weta(t),*} \geq \delta_\rho \quad \hbox{for all}\ t\in [0,1].
    $$
By our choice of $\varphi$, $\psi$, we have
    $$  \eqalign{
    {d\over dt}J(\weta(t)) &= \Diff J(\weta(t)) {d\weta\over dt}(t) 
     = -\Diff J(\weta(t)) {\calV(\weta(t))\over \norm{\calV(\weta(t))}_{\weta(t)}} \cr
    &\leq  -\half \norm{\Diff J(\weta(t))}_{\weta(t),*} 
    \leq  -\half \delta_\rho.               \cr}
    $$
Thus we have
    $$  J(\weta(1)) = J(\weta(0)) + \int_0^1 {d\over dt}J(\weta(t))\, dt
        \leq J(\weta(0)) -{\delta_\rho\over 2} \leq b+\epsilon -{\delta_\rho\over 2}.
    $$
If Case 2 takes a place, we can find an interval $[\alpha,\beta]\subset [0,1]$ such that
    $$  \eqalign{
        &\weta(\alpha) \in \partial \wN_\rho(\wKb), \quad
        \weta(\beta) \in \partial \wN_{{2\over 3}\rho}(\wKb), \cr
        &\weta(t) \in \wN_\rho(\wKb)\setminus \wN_{{2\over 3}\rho}(\wKb) \quad
            \hbox{for all}\ t\in [\alpha,\beta).   \cr}
    $$
By \ref[4.21],
    $$  \beta-\alpha \geq \int_\alpha^\beta \norm{{d\weta\over dt}(t)}_{\weta(t)}\, dt 
        \geq \distM(\weta(\alpha),\weta(\beta)) \geq {1\over 3}\rho. 
    $$
Thus,
    $$  \eqalign{
        J(\weta(1)) & \leq J(\weta(\beta)) = J(\weta(\alpha)) + \int_\alpha^\beta {d\over dt}J(\weta(t))\, dt \cr
        & \leq J(\weta(0)) + \int_\alpha^\beta {d\over dt}J(\weta(t))\, dt 
          \leq J(\weta(0)) + \int_\alpha^\beta -{\delta_\rho\over 2} \, dt \cr
        & \leq J(\weta(0)) -{\delta_\rho\over 2}(\beta-\alpha) 
          \leq b+\epsilon -{\delta_\rho \rho\over 6}.\cr}
    $$
Choosing $\epsilon < \min\{{\overline\epsilon\over 2}, {\delta_\rho\over 4}, {1\over 12}\delta_\rho \rho\}$,
we have $J(\weta(1)) \leq b-\epsilon$ in both cases.  This is a contradiction and we have 
\ref[4.20].
                                                                            \QED

\medskip

In the following section, we can construct a deformation flow for $I(\lambda,u)$ using\m
$\weta(t,\tlu)$.

\medskip


\BSScap{\label[Subsection:4.d].  Deformation flow for $I(\lambda,u)$}{\raw[4.d].  Deformation flow for I}
In this section, we construct a deformation flow for $I(\lambda,u)$ and give a proof to
our Proposition \ref[Proposition:3.1].

We use the following maps:
    $$  \eqalign{
    &\pi:\, M\to \RE;\ (\tlu(x)) \mapsto (\lambda, u(x/e^\theta)), \cr
    &\iota:\, \RE\to M;\ (\lambda,u(x))\mapsto (0,\lambda,u(x)) \cr}
    $$
and we construct a deformation flow $\eta(t,\lambda,u):\, [0,1]\times\RE\to \RE$ as 
a composition $\pi\circ\weta(t,\cdot)\circ\iota$;
    $$  \eta(t,\lambda,u)=\pi(\weta(t,\iota(\lambda,u)))
        = \pi(\weta(t,0,\lambda,u)).                    \eqno\label[4.22]
    $$
As fundamental properties of $\pi$ and $\iota$, we have
    $$  \eqalign{
    &\pi(\iota(\lambda,u)) = (\lambda,u) \quad \hbox{for all}\ (\lambda,u)\in\RE,     \cr
    &\iota(\pi(\tlu)) = (0,\lambda, u(x/e^\theta)) \quad \hbox{for all}\ (\tlu)\in M, \cr
    &J(\tlu)=I(\pi(\tlu)) \quad \hbox{for all}\ (\tlu)\in M. \cr}
    $$
Clearly $\pi(\wKb)=K_b$.  The following lemma gives us a relation between $\pi(\wN_\rho(\wKb))$
and $N_\rho(K_b)$.

\medskip

\proclaim Lemma \label[Lemma:4.9].
For any $\rho>0$ there exists $R(\rho)>0$ such that
    $$  \eqalignno{
        &\pi(\wN_\rho(\wKb)) \subset N_{R(\rho)}(K_b),      &\label[4.23]\cr
        &\iota\Bigl(((\RE)\setminus N_{R(\rho)}(K_b))\Bigr) \subset M\setminus \wN_\rho(\wKb).  &\label[4.24]\cr}
    $$
Moreover
    $$  R(\rho)\to 0 \quad \hbox{as}\ \rho\to 0.            \eqno\label[4.25]
    $$

\claim Proof.
For $\rho>0$, suppose that $(\lambda_0,u_0)\in\RE$ satisfies $\distM((0,\lambda_0,u_0),\wKb)\leq\rho$.
First we show
    $$  \distM((\lambda_0,u_0), K_b) \leq e^{N\rho\over 2} \rho
            +\sup\{ \norm{\omega(e^\alpha x)-\omega(x)}_{H^1};\,
            \abs\alpha\leq \rho,\, \omega\in P_2(K_b)\}.               \eqno\label[4.26]
    $$
In fact, for any $\epsilon>0$ there exists $\sigma(t)=(\theta(t),\lambda(t), u(t))\in 
C^1([0,1],M)$ such that $\sigma(0)=(0,\lambda_0,u_0)$, $\sigma(1)\in \wKb$ and
    $$  \int_0^1 \norm{\dot\sigma(t)}_{\sigma(t)} \, dt \leq \rho+\epsilon.
    $$
In particular, since $\theta(0)=0$, for any $t\in [0,1]$
    $$  \abs{\theta(t)} \leq \int_0^1\abs{\dot\theta(t)}\, dt 
        \leq \int_0^1 \norm{\dot\sigma(t)}_{\sigma(t)}\, dt 
        \leq \rho+\epsilon.
    $$
Thus
    $$  \eqalign{
    &\norm{(\lambda(0),u(0))-(\lambda(1),u(1))}_{\RE}
    \leq \int_0^1 \left(\abs{\dot\lambda(t)}^2 + \norm{\dot u(t)}_{H^1}^2\right)^{1/2}\, dt \cr
    \leq& e^{N(\rho+\epsilon)\over 2}\int_0^1 
        \left(\abs{\dot\theta(t)}^2 + \abs{\dot\lambda(t)}^2 
            + e^{(N-2)\theta(t)}\norm{\nabla \dot u(t)}_2^2 
            + e^{N\theta(t)}\norm{\dot u(t)}_2^2 \right)^{1/2} \, dt \cr
    =& e^{N(\rho+\epsilon)\over 2}\int_0^1 
            \norm{(\dot\theta(t),\dot\lambda(t),\dot u(t))}_{\sigma(t)}\, dt
    \leq e^{N(\rho+\epsilon)\over 2} (\rho+\epsilon).\cr}
    $$
On the other hand, since $(\theta(1),\lambda(1),u(1))\in \wKb$, we have
$(\lambda(1),u(1)(x/e^{\theta(1)}))\in K_b$, i.e., $u(1)(x/e^{\theta(1)})\in P_2(K_b)$.
Thus
    $$  \eqalign{
        &\dist((\lambda_0,u_0),K_b) 
        \leq \norm{(\lambda(0),u(0))-(\lambda(1),u(1)(x/e^{\theta(1)}))}_{\RE} \cr
        \leq& \norm{(\lambda(0),u(0))-(\lambda(1),u(1))}_{\RE} \cr
            &+\norm{(\lambda(1),u(1)(x))-(\lambda(1),u(1)(x/e^{\theta(1)}))}_{\RE} \cr
        \leq& \norm{(\lambda(0),u(0))-(\lambda(1),u(1))}_{\RE}
            +\norm{u(1)(x)-u(1)(x/e^{\theta(1)})}_{H^1} \cr
        \leq& e^{N(\rho+\epsilon)\over 2}(\rho+\epsilon)
            +\sup\{\norm{\omega(e^\alpha x)-\omega(x)}_{H^1};\, 
            \abs{\alpha}\leq \rho+\epsilon, \,  \omega\in P_2(K_b)\}.  \cr}
    $$
Since $\epsilon>0$ is arbitrary, we have \ref[4.26].

We set 
    $$  R(\rho)= e^{N\rho\over 2}\rho +\sup\{\norm{\omega(e^\alpha x)-\omega(x)}_{H^1};\, 
            \abs{\alpha}\leq \rho, \,   \omega\in P_2(K_b)\}.
    $$
Then
    $$  \distM((0,\lambda_0,u_0),\wKb) \leq\rho \quad \hbox{implies}\quad
        \dist((\lambda_0,u_0),K_b) \leq R(\rho).            \eqno\label[4.27]
    $$
Since $P_2(K_b)$ is compact in $\E$, we have
    $$  \sup\{\norm{\omega(e^\alpha x)-\omega(x)}_{H^1};\, 
            \abs{\alpha}\leq \rho, \,   \omega\in P_2(K_b)\} \to 0
        \quad \hbox{as}\ \rho\to 0,
    $$
which implies \ref[4.25].  
Noting $\distM((\tlu),\wKb) = \distM((0,\lambda, u(x/e^\theta)), \wKb)$, 
\ref[4.27] implies \ref[4.23] and \ref[4.24].  \QED

\medskip

Now we can give a proof of Proposition \ref[Proposition:3.1].

\medskip

\claim Proof of Proposition \ref[Proposition:3.1].
Let $\calO$ be a given neighborhood of $K_b$ and let $\overline\epsilon>0$ be a given
positive number.

We take small $\rho>0$ such that $N_{R(\rho)}(K_b) \subset\calO$.  By Proposition \ref[Proposition:4.8],
there exist $\epsilon\in (0,\overline\epsilon)$ and $\weta:\, [0,1]\times M\to M$ such that
\maru 1--\maru 6 in Proposition \ref[Proposition:4.8] hold.
We define $\eta(t,\lambda,u):\, [0,1]\times\RE\to \RE$ by \ref[4.22].

We can check that $\eta(t,\lambda,u)$ satisfies the properties \maru 1--\maru 6 of 
Proposition \ref[Proposition:3.1].  Here we just prove
    $$  \eta(1,[I\leq b+\epsilon]\setminus\calO) \subset [I\leq b-\epsilon].    \eqno\label[4.28]
    $$
Since $[I\leq b+\epsilon]\setminus\calO\subset [I\leq b+\epsilon]\setminus N_{R(\rho)}(K_b)$,
we have from \ref[4.24]
    $$  \iota([I\leq b+\epsilon]\setminus\calO)\subset [J\leq b+\epsilon]_M\setminus \wN_\rho(\wKb).
                                                                \eqno\label[4.29]
    $$
By \maru 4 of Proposition \ref[Proposition:4.8],
    $$  \weta(1,[J\leq b+\epsilon]_M\setminus \wN_\rho(\wKb)) \subset [J\leq b-\epsilon]_M.
                                                                \eqno\label[4.30]
    $$
By the definition of $\pi$ and \ref[4.6],
    $$  \pi([J\leq b-\epsilon]_M) \subset [I\leq b-\epsilon].   \eqno\label[4.31]
    $$
Combining \ref[4.29]--\ref[4.31], we have \ref[4.28].  \QED

\medskip

\claim Remark \label[Remark:4.10].
By our construction, $t\mapsto \weta(t,\theta,\lambda,u);\, [0,1]\to \RRE$ is of class $C^1$.
However, $t\mapsto u(x/e^t);\, \R\to \E$ is continuous but not of class $C^1$ for
$u\in \E\setminus H^2(\R^N)$ and thus $t\mapsto \eta(t,\lambda,u)=\pi(\weta(t,0,\lambda,u));
\, [0,1]\to\RE$ is continuous but not of class $C^1$.

\medskip

\BS{\label[Section:5]. Minimizing problem}
In this section we assume \cond[g:1]--\cond[g:4] (without \cond[g:5]) and we deal
with Theorem \ref[Theorem:0.1].
Under the condition $\calI_m<0$, the existence of a solution is shown by Shibata [\cite[S1]],
that is, he showed that $\calI_m$ is achieved by a solution of $(*)_m$.  First we give 
an approach using our functional $I(\lambda,u)$.

\BSS{\label[Subsection:5.a].  Mountain pass approach}
Under the condition \cond[g:1]--\cond[g:4], as in Sections \ref[Section:1]--\ref[Section:2], we define 
$\lambda_0\in (-\infty,\infty]$ by \ref[1.8].

\smallskip

\item{(i)} For $\lambda<\lambda_0$, $u\mapsto \whI(\lambda,u)$ has the mountain pass geometry.
\item{(ii)} When $\lambda_0<\infty$, $\whI(\lambda,u)\geq 0$ for all $\lambda\geq\lambda_0$ and
$u\in\E$.

\smallskip

\noindent
We set for $\lambda<\lambda_0$
    $$  \eqalignno{
    &\whGamma_{mp}(\lambda) = \{ \zeta(\tau)\in C([0,1],\E);\, \zeta(0)=0,\, 
        \whI(\lambda,\zeta(1))<0\},         \cr
    &a_{mp}(\lambda) = \inf_{\zeta\in\whGamma_{mp}(\lambda)}\max_{\tau\in [0,1]}\whI(\lambda,\zeta(\tau)).
                                            &\label[5.1]\cr}
    $$
We note that if \cond[g:5] holds, $a_{mp}(\lambda)$ coincides with $a_1(\lambda)$ defined in 
\ref[1.9]--\ref[1.10].
By the result of [\cite[HIT]], we see that $a_{mp}(\lambda)$ is attained by a critical point of
$u\mapsto \whI(\lambda,u)$.
This fact can also be shown via our new deformation argument.  See Section \ref[Section:6].

We set
    $$  m_0 = 2\inf_{\lambda\in (-\infty,\lambda_0)}{ a_{mp}(\lambda)\over e^\lambda}.
                                                            \eqno\label[5.2]
    $$
As in Sections \ref[Section:1]--\ref[Section:4], we can show

\medskip

\proclaim Theorem \label[Theorem:5.1].
Assume \cond[g:1]--\cond[g:4].  Suppose $m>m_0$.  Then $(*)_m$ has at least one solution
$(\lambdamp,\ump)$, which is characterized by the following minimax method;
    $$  I(\lambdamp,\ump) = b_{mp}<0,
    $$
where 
    $$  \eqalignno{
        &b_{mp} = \inf_{\gamma\in \Gamma_{mp}}\max_{\tau\in[0,1]} I(\gamma(\tau)), \cr
        &\Gamma_{mp} = \{\gamma(\tau)\in C([0,1],\RE);\, 
            \gamma(0)\in [\lambda_m,\infty)\times\{0\},\, I(\gamma(0)) \leq B_m-1, \cr
        &\qquad\qquad\qquad\qquad
            \gamma(1)\in (\RE)\setminus\Omega_m,\, I(\gamma(1)) \leq B_m-1\}.       \cr}
    $$
Here $\lambda_m\in\R$, $\Omega_m\subset [\lambda_m,\infty)\times\E$ and 
$B_m = \inf_{(\lambda,u)\in \partial\Omega_m} I(\lambda,u) > -\infty$ are chosen 
as in Section \ref[Section:2]. \QED

\medskip

\noindent
As a corollary, we have

\medskip

\proclaim Corollary \label[Corollary:5.2].
Assume \cond[g:1]--\cond[g:4] and suppose $m>m_0$.  Then 
    $$  \calI_m < 0.
    $$

\claim Proof.  The critical point $(\lambdamp,\ump)$ obtained in Theorem \ref[Theorem:5.1] satisfies 
    $$  \norm{\ump}_2^2 =m \quad \hbox{and}\quad \calF(\ump) = I(\lambdamp,\ump) = b_{mp}<0.
    $$
Thus $\calI_m = \inf_{\norm u_2^2 =m} \calF(u) \leq \calF(\ump) <0$.  \QED

\medskip

We also have

\medskip

\proclaim Theorem \label[Theorem:5.3].
Under the assumption of Theorem \ref[Theorem:5.1], there exists $\gamma_0\in\Gamma_{mp}$ such that
    $$  b_{mp} = \max_{\tau\in [0,1]} I(\gamma_0(\tau)).
    $$

\claim Proof.
Let $(\lambdamp,\ump)$ be the critical point corresponding to $b_{mp}$.
In Jeanjean-Tanaka [\cite[JT]], we find a path $\zeta_0(\tau)\in \whGamma_{mp}(\lambdamp)$
such that
    $$  \ump \in \zeta_0([0,1]) \quad \hbox{and} \quad
        b_{mp} = \whI(\lambdamp,\ump) 
            = \max_{\tau\in [0,1]} \whI(\lambdamp,\gamma_0(\tau)). 
    $$
As in the proof of Lemma \ref[Lemma:2.3], we may assume $\whI(\lambdamp,\zeta_0(1))\leq B_m-1$.
Joining paths
    $$  [0,1]\to \RE;\, \tau\mapsto (\lambdamp\tau + L(1-\tau), 0)
    $$
and 
    $$  [0,1]\to \RE;\, \tau\mapsto (\lambdamp, \zeta_0(\tau)),
    $$
we find the desired path $\gamma_0\in \Gamma_{mp}$.  \QED

\medskip


\BSScap{\label[Subsection:5.b]. Mountain pass characterization of $\calI_m$}{\raw[5.b]. Mountain pass characterization of Im}
Next we consider the problem $(*)_m$ under the conditions \cond[g:1]--\cond[g:4] and $\calI_m<0$.

Shibata [\cite[S1]] showed the following

\medskip

\proclaim Theorem \label[Theorem:5.4] {\rm ([\cite[S1]])}.
There exists $m_S\in [0,\infty)$ such that
\item{(i)} $\calI_m=0$ for $m\in (0,m_S]$,\m
$\calI_m<0$ for $m\in (m_S,\infty)$.
\item{(ii)} If $\calI_m<0$, $\calI_m$ is attained and the minimizer is a solution
of $(*)_m$.   \QED

\medskip

\noindent
In what follows, we will show that $m_0$ given in \ref[5.2] coincides with $m_S$ and
$\calI_m=b_{mp}$.  Precisely,

\smallskip

\item{(i)} $m > m_0$ if and only if $\calI_m <0$.
\item{(ii)} For $m>m_0$,  $\calI_m = b_{mp}$.

\smallskip

\noindent
First we show the minimizer of $\calI_m$ satisfies the following properties.

\medskip

\proclaim Lemma \label[Lemma:5.5].
Suppose $\calI_m<0$ and let $(\mumin,\umin)$ be the corresponding minimizer of $\calI_m$, i.e.,
$\calF(\umin)=\calI_m$, $\norm{\umin}_2^2=m$.  Then
\item{(i)} $\mumin>0$.
\item{(ii)} ${N-2\over 2}\norm{\nabla\umin}_2^2 +N\left( {\mumin\over 2}\norm{\umin}_2^2
-\intRN G(\umin)\right)=0$.

\medskip

\claim Proof.
First we show (ii).  We set $u_{*\theta}(x)=\theta^{N/2}\umin(\theta x)$ for $\theta>0$.
Since $\umin$ is a minimizer of $\calF(u)$ under the constraint $\norm u_2^2=m$ and 
$\norm{u_{*\theta}}_2^2 =m$ for all $\theta>0$, we have
    $$  {d\over d\theta}\Bigr|_{\theta=1} \calF(u_{*\theta}) = 0,
    $$
that is,
    $$  \norm{\nabla\umin}_2^2 + N\intRN G(\umin) -{N\over 2}\intRN g(\umin)\umin=0.    \eqno\label[5.3]
    $$
Since $(\mumin,\umin)$ solves $(*)_m$, we also have
    $$  \norm{\nabla\umin}_2^2 +\mumin\norm{\umin}_2^2 = \intRN g(\umin)\umin.      \eqno\label[5.4]
    $$
(ii) follows from \ref[5.3] and \ref[5.4].

Next we show (i).  By (ii), we have
    $$  {\mumin N\over 2}m = {\mumin N\over 2}\norm{\umin}_2^2  = -N\calF(\umin) +\norm{\nabla\umin}_2^2 
        \geq -N\calI_m >0. 
    $$
Thus we have $\mumin>0$.                            \QED
     
\medskip

\noindent
By Lemma \ref[Lemma:5.5], setting $\lambdas=\log\mumin$, $(\lambdas,\umin)$ is a critical point of $I(\lambda,u)$
with
    $$  I(\lambdas,\umin)=\calI_m \quad \hbox{and}\quad P(\lambdas,\umin)=0.
    $$
Next we show

\medskip

\proclaim Proposition \label[Proposition:5.6].
Suppose $\calI_m<0$ and let $(\lambdas,\umin)$ be a critical point corresponding to $\calI_m$.
Then we have
\item{(i)} $u\mapsto \whI(\lambdas,u)$ has the mountain pass geometry, that is, $\lambdas<\lambda_0$.
\item{(ii)} $\whI(\lambdas,\umin)\geq a_{mp}(\lambdas)$.
\item{(iii)} $m>m_0$, where $m_0$ is given in \ref[5.2].

\medskip

\claim Proof.
(i) It suffices to show $\whI(\lambdas,u)<0$ for some $u\in \E$.  We set 
    $$  \wG(\xi) = G(\xi) -{e^{\lambdas}\over 2}\xi^2 \quad \hbox{for}\ \xi\in\R.
    $$
Then we have for some $v\in\E$
    $$  \intRN \wG(v) >0.                           \eqno\label[5.5]
    $$
In fact, when $N\geq 3$, it follows from $P(\lambdas,\umin)=0$ that \ref[5.5] holds with
$v=\umin$.  \m
When $N=2$, we have by $P(\lambdas,\umin)=0$
    $$  \intRN \wG(\umin)=0.
    $$
We also have from $(*)_m$
    $$  {d\over ds}\Bigr|_{s=1} \intRN\wG(s\umin) = \intRN g(\umin)\umin -e^{\lambdas}\norm{\umin}_2^2
        = \norm{\nabla\umin}_2^2 >0.
    $$
Thus \ref[5.5] holds with $v=s\umin$ for $s>1$ closed to $1$. Since 
    $$  \whI(\lambdas, v(x/\theta)) = \half \theta^{N-2}\norm{\nabla v}_2^2 -\theta^N\intRN \wG(v)
        < 0 \quad \hbox{for large}\ \theta\gg 1,
    $$
(i) holds. \m
(ii) By the result of [\cite[JT]], the mountain pass minimax value $a_{mp}(\lambdas)$
gives the least energy level for $\whI(\lambdas,u)$.  Thus 
$\whI(\lambdas,\umin) \geq a_{mp}(\lambdas)$. \m
(iii)  (ii) implies
    $$  \eqalign{
        {e^{\lambdas}\over 2}m &= {e^{\lambdas}\over 2}\norm{\umin}_2^2 
            =\whI(\lambdas,\umin) -\calF(\umin) \cr
        &\geq a_{mp}(\lambdas)-\calI_m > a_{mp}(\lambdas).  \cr}
    $$
Thus
    $$  m > 2{a_{mp}(\lambdas)\over e^{\lambdas}} \geq m_0.     \QED
    $$

\medskip

\proclaim Proposition \label[Proposition:5.7].
Suppose $\calI_m<0$.  Then  $\calI_m=b_{mp}$.

\medskip

\claim Proof.
As in Lemma \ref[Lemma:2.3], we can show
    $$  b_{mp} \leq a_{mp}(\lambda) -{e^\lambda\over 2}m 
       \quad \hbox{for all}\ \lambda\in (-\infty,\lambda_0).
    $$
Thus, by (ii) of Proposition \ref[Proposition:5.6], for a critical point $(\lambdas,\umin)$ corresponding to 
$\calI_m$,
    $$  \calI_m =I(\lambdas,\umin) = \whI(\lambdas,\umin) - {e^{\lambdas}\over 2}m
        \geq a_{mp}(\lambdas) - {e^{\lambdas}\over 2}m \geq b_{mp}.
    $$
On the other hand, it follows from (iii) of Proposition \ref[Proposition:5.6] that $m>m_0$ and $b_{mp}$ is attained 
by a critical point $(\lambdamp,\ump)\in\RE$.  Thus
    $$  \norm{\ump}_2^2 = m, \quad \calF(\ump) = I(\lambdamp,\ump) = b_{mp}.
    $$
Thus
    $$  \calI_m = \inf_{\norm u_2^2=m}\calF(u) \leq \calF(\ump) = b_{mp}.
    $$
Therefore we have $\calI_m=b_{mp}$.                         \QED

\medskip

\noindent
We also have

\medskip

\proclaim Corollary \label[Corollary:5.8].
$\calI_m<0$ if and only if $m>m_0$.

\medskip

\claim Proof.
``if'' part follows from Theorem \ref[Theorem:5.1] and
``only if'' part follows from Proposition \ref[Proposition:5.6].  \QED

\medskip

\claim End of the proof of Theorem \ref[Theorem:0.1].
Theorem \ref[Theorem:0.1] follows from Theorem \ref[Theorem:5.1], Propositions \ref[Proposition:5.6], \ref[Proposition:5.7] and Corollary \ref[Corollary:5.8].  \QED

\medskip

\BS{\label[Section:6]. Deformation lemma for scalar field equations}
In this section we study the following nonlinear scalar field equations:
    $$  \left\{ \eqalign{
        -&\Delta u + \mu u = g(u) \quad \hbox{in}\ \R^N, \cr
        &u\in H^1(\R^N),  \cr}
        \right.                 \eqno\label[6.1]
    $$
where $N\geq 2$, $\mu>0$ and $g(\xi)\in C(\R,\R)$ satisfies \cond[g:1], \cond[g:2], 
\cond[g:3] with $p={N+2\over N-2}$ ($N\geq 3$), $p\in (1,\infty)$ ($N=2$).
Solutions of \ref[6.1] are characterized as critical points of the following functional:
    $$  I(u) = \half\norm{\nabla u}_2^2 +{\mu\over 2}\norm u_2^2 -\intRN G(u) \in C^1(\E,\R).
    $$
Here we use notation different from previous sections.  We also write
    $$  P(u) = {N-2\over 2}\norm{\nabla u}_2^2 +N\left({\mu\over 2}\norm u_2^2 -\intRN G(u)\right).
    $$
In this section we give a new deformation result for \ref[6.1] using ideas in 
Sections \ref[Section:3]--\ref[Section:4].

A key of our argument is the following

\medskip

\proclaim Proposition \label[Proposition:6.1].
For any $b\in\R$, $I(u)$ satisfies the following $(PSP')_b$:
\itemitem{$(PSP')_b$} If a sequence $(u_n)_{n=1}^\infty\subset \E$ satisfies as
$n\to\infty$
    $$  \eqalignno{
        &I(u_n) \to b,                          &\label[6.2]\cr
        &\partial_u I(u_n) \to 0 \quad \hbox{strongly in}\ (\E)^*,          &\label[6.3]\cr
        &P(u_n) \to 0,                          &\label[6.4]\cr}
    $$
then $(u_n)_{n=1}^\infty$ has a strongly convergent subsequence in $\E$.

\claim Proof.
First we note by \cond[g:2], \cond[g:3] with $p={N+2\over N-2}$ ($N\geq 3$), $p\in (0,\infty)$ ($N=2$) that
$u_m\wlimit u_0$ weakly in $\E$ implies for any $\varphi\in\E$
    $$  \intRN g(u_n)\varphi \to \intRN g(u_0)\varphi, \quad \intRN g(u_n)u_n \to \intRN g(u_0)u_0
                                            \eqno\label[6.5]
    $$
Proof consists of several steps.  Here we follow essentially the argument 
in [\cite[HIT]] (Propositions \raw[5.1] and \raw[5.3]).

\smallskip

\noindent
{\sl Step 1: $\norm{\nabla u_n}_2$ is bounded as $n\to \infty$.
}

\smallskip

\noindent
Since $\norm{\nabla u_n}_2^2 =NI(u_n)-P(u_n)$, Step 1 follows from \ref[6.2] and \ref[6.4].

\smallskip

\noindent
From now on we prove that $\norm{u_n}_2$ is bounded as $n\to \infty$.  We argue indirectly
and we assume
    $$  t_n = \norm{u_n}_2^{-2/N} \to 0 \quad \hbox{as}\ n\to \infty.
    $$
We set $v_n(x)=u_n(x/t_n)$.  Since 
    $$  \norm{v_n}_2^2 = 1 \quad \hbox{and}\quad \norm{\nabla v_n}_2^2 = t_n^{N-2}\norm{\nabla u_n}_2^2,
                                            \eqno\label[6.6]
    $$
$(v_n)_{n=1}^\infty$ is bounded in $\E$.  Thus we may assume after extracting a subsequence that
    $$  v_n \wlimit v_0 \quad \hbox{weakly in}\  \E.
    $$

\noindent
{\sl Step 2: $v_0=0$.
}

\smallskip

\noindent
Denoting $\epsilon_n\equiv \norm{\partial_u I(u_n)}_{(\E)^*}\to 0$, we have
    $$  \pabs{ (\nabla u_n,\nabla \zeta)_2 +\mu(u_n, \zeta)_2 -\intRN g(u_n)\zeta} 
        \leq \epsilon_n \norm{\zeta}_{H^1} \quad \hbox{for any}\ \zeta\in \E.
    $$
Setting $u_n(x)=v_n(t_nx)$, $\zeta(x)=\varphi(t_nx)$, where $\varphi\in\E$,
    $$  \eqalign{
        &\pabs{ t_n^{-(N-2)} (\nabla v_n,\nabla \varphi)_2 +\mu t_n^{-N}(v_n, \varphi)_2 - t_n^{-N}\intRN g(v_n)\varphi} \cr
        &\leq \epsilon_n \left(t_n^{-(N-2)}\norm{\nabla \varphi}_2^2 + t_n^{-N}\norm{\varphi}_2^2\right)^{1/2}. \cr}
    $$
Thus
    $$  \pabs{ t_n^2 (\nabla v_n,\nabla \varphi)_2 +\mu (v_n, \varphi)_2 - \intRN g(v_n)\varphi} 
        \leq \epsilon_n t_n^{N/2} \left(t_n^2\norm{\nabla \varphi}_2^2 + \norm{\varphi}_2^2\right)^{1/2}, \eqno\label[6.7]
    $$
from which we have
    $$  \intRN (\mu v_0-g(v_0))\varphi =0 \quad \hbox{for any}\ \varphi\in \E.
    $$
Thus $\mu v_0-g(v_0) =0$.  Since $\xi=0$ is an isolated solution of $\mu\xi-g(\xi)=0$ by \cond[g:2], we have
$v_0(x)\equiv 0$.

\smallskip

\noindent
{\sl Step 3: $\norm{u_n}_2$ is bounded as $n\to \infty$.
}

\smallskip

\noindent
Setting $\varphi=v_n$ in \ref[6.7],
    $$  \pabs{ t_n^2 \norm{\nabla v_n}_2^2 +\mu \norm{v_n}_2^2 - \intRN g(v_n)v_n} 
        \leq \epsilon_n t_n^{N/2} \left(t_n^2\norm{\nabla v_n}_2^2 + \norm{v_n}_2^2\right)^{1/2}. 
    $$
Thus, by \ref[6.5], $\norm{v_n}_2\to 0$ as $n\to \infty$, which contradicts with \ref[6.6].
Thus $(u_n)_{n=1}^\infty$ is bounded in $\E$.

\smallskip

\noindent
{\sl Step 4: Conclusion.
}

\smallskip

\noindent
By Step 1 and Step 3, $(u_n)_{n=1}^\infty$ is bounded in $\E$.  
After extracting a subsequence, we may assume that $u_n\wlimit u_0$ weakly in $\E$ for some $u_0$.
Since $\partial_u I(u_n)u_n\to 0$, $\partial_u I(u_n)u_0\to 0$, we deduce from \ref[6.5] that
    $$  \lim_{n\to\infty} (\norm{\nabla u_n}_2^2 +\mu\norm{u_n}_2^2) = \norm{\nabla u_0}_2^2 +\mu\norm{u_0}_2^2.
    $$
Thus $u_n\to u_0$ strongly in $\E$.  \QED

\medskip

Arguing as in Sections \ref[Section:3]--\ref[Section:4], we have

\medskip

\proclaim Proposition \label[Proposition:6.2].
Under the assumption of Proposition \ref[Proposition:6.1], for any $b\in\R$ we have
\item{(i)} $K_b = \{ u\in\E;\, I(u)=b, \, \partial_u I(u)=0, \, P(u)=0\}$ 
is compact in $\E$.
\item{(ii)} For any open neighborhood $\calO$ of $K_b$ and $\overline\epsilon>0$ there exist
$\epsilon\in (0,\overline\epsilon)$ and a continuous map $\eta(t,u):\, [0,1]\times\E
\to\E$ such that
\itemitem{\maru 1} $\eta(0,u)=u$ for all $u\in\E$.
\itemitem{\maru 2} $\eta(t,u)=u$ if  $u\in [I\leq b-\overline\epsilon]$.
\itemitem{\maru 3} $I(\eta(t,u))\leq I(u)$ for all $(t,u)\in [0,1]\times\E$.
\itemitem{\maru 4} $\eta(1,[I\leq b+\epsilon]\setminus\calO)\subset [I\leq b-\epsilon]$, \m
        $\eta(1,[I\leq b+\epsilon])\subset [I\leq b-\epsilon]\cup\calO$.
\itemitem{\maru 5} If $K_b=\emptyset$, $\eta(1,[I\leq b+\epsilon])\subset [I\leq b-\epsilon]$.
\item{}
Here we use notation: $[I\leq c] = \{ u\in\E;\, I(u)\leq c\}$ for $c\in\R$.

\medskip

\noindent
Using Proposition \ref[Proposition:6.2], we can show that $a_{mp}(\lambda)$ given in \ref[5.1] is a critical 
value of $u\mapsto \whI(\lambda,u)$.

\medskip

\BS{\labeldef[Section:A]. Appendix: Proof of (iv) of Lemma \ref[Lemma:2.4]}
In this appendix we give a proof to (iv) of Lemma \ref[Lemma:2.4].

\medskip

\claim Proof of (iv) of Lemma \ref[Lemma:2.4].
Suppose that a closed set $Z$ is invariant under $\Z_2$-action \ref[2.6]
and satisfies $0\not\in \overline{P_2(Z)}$.  Then $\overline{P_2(Z)}\subset \E$
is symmetric with respect to $0$ and $\g(\overline{P_2(Z)})$ is well-defined.

For $A=\gamma(\overline{D_{j+\ell}\setminus Y})$, $\gamma\in\Gamma_{j+\ell}$,
$\g(Y)\leq \ell$, we have
    $$  \overline{A\setminus Z} 
    =\gamma(\overline{D_{j+\ell}\setminus (Y\cup\gamma^{-1}(Z))}). \eqno\label[A.1]
    $$
In fact,
    $$  \gamma(D_{j+\ell}\setminus (Y\cup\gamma^{-1}(Z))) 
        = \gamma(D_{j+\ell}\setminus Y)\setminus Z
        \subset \gamma(\overline{D_{j+\ell}\setminus Y})\setminus Z = A\setminus Z.
                                \eqno\label[A.2]
    $$
Conversely, since $\overline{B}\setminus C\subset\overline{B\setminus C}$ for a set $B$
and a closed set $C$, we have
    $$  A\setminus Z =\gamma(\overline{D_{j+\ell}\setminus Y})\setminus Z
        = \gamma(\overline{D_{j+\ell}\setminus Y}\setminus \gamma^{-1}(Z))
        \subset \gamma(\overline{D_{j+\ell}\setminus (Y\cup\gamma^{-1}(Z))}).
                                            \eqno\label[A.3]
    $$
Thus \ref[A.1] follows from \ref[A.2] and \ref[A.3].
Since $P_2\circ \gamma:\, \gamma^{-1}(Z)\to \overline{P_2(Z)}$ is an odd map,
    $$  \g(\gamma^{-1}(Z)) \leq \g(\overline{P_2(Z)}) =i.
    $$
Thus,
    $$  \eqalign{
        \g(Y\cup\gamma^{-1}(Z)) &\leq \g(Y) +\g(\gamma^{-1}(Z)) 
            \leq \g(Y) +\g(\overline{P_2(Y)}) \cr
            &\leq \ell+i. \cr}
    $$
Therefore, by \ref[A.1] we have $\overline{A\setminus Z}\in \Lambda_{j-i}$.  \QED

\medskip

\centerline{\bf Acknowledgments}

\smallskip

\noindent
The authors would like to thank Professor Tohru Ozawa for helpful discussions and
comments.

\medskip

\bibliography

\medskip

\bibitem[AdAP] 
A. Azzollini, P. d'Avenia, A. Pomponio, 
Multiple critical points for a class of nonlinear functionals, 
Ann. Mat. Pura Appl. (4) 190 (2011), no. 3, 507--523.

\bibitem[BV]
T. Bartsch, S. de Valeriola, 
Normalized solutions of nonlinear Schr\"odinger equations, 
Arch. Math. (Basel) 100 (2013), no. 1, 75--83.

\bibitem[BGK]
H. Berestycki, T. Gallou\"et, O. Kavian, 
\'Equations de champs scalaires euclidiens non lin\'eaires dans le plan, 
C. R. Acad. Sci. Paris Ser. I Math. 297 (1983), no. 5, 307--310 and
Publications du Laboratoire d'Analyse Numerique, Universite de Paris VI, (1984).

\bibitem[BL1]
H. Berestycki, P.-L. Lions, Nonlinear scalar field equations. I. Existence of a ground state, 
Arch. Rational Mech. Anal. 82 (1983), no. 4, 313--345.

\bibitem[BL2]
H. Berestycki, P.-L. Lions, 
Nonlinear scalar field equations. II. Existence of infinitely many solutions, 
Arch. Rational Mech. Anal. 82 (1983), no. 4, 347--375.

\bibitem[BT]
J. Byeon, K. Tanaka, 
Nonlinear elliptic equations in strip-like domains, 
Advanced Nonlinear Studies 12 (2012), 749--765.

\bibitem[CL]
T. Cazenave, P.-L. Lions, 
Orbital stability of standing waves for some nonlinear Schr\"odinger equations, 
Comm. Math. Phys. 85 (1982), no. 4, 549--561.

\bibitem[CT]
C.-N. Chen, K. Tanaka, 
A variational approach for standing waves of FitzHugh--Nagumo type systems, 
J. Differential Equations 257 (2014), no. 1, 109--144.

\bibitem[HIT]
J. Hirata, N. Ikoma, K. Tanaka, 
Nonlinear scalar field equations in {$\R^N$}: mountain pass and symmetric mountain 
pass approaches, 
Topol. Methods Nonlinear Anal. 35 (2010), no. 2, 253--276.

\bibitem[J]
L. Jeanjean, 
Existence of solutions with prescribed norm for semilinear elliptic equations, 
Nonlinear Anal. 28 (1997), no. 10, 1633--1659.

\bibitem[JT]
L. Jeanjean, K. Tanaka, 
A remark on least energy solutions in $\R^N$, 
Proc. Amer. Math. Soc. 131 (2002), Number 8, Pages 2399--2408.

\bibitem[MVS]
V. Moroz, J. Van Schaftingen, 
Existence of groundstates for a class of nonlinear Choquard equations, 
Trans. Amer. Math. Soc. 367 (2015), no. 9, 6557--6579.

\bibitem[R]
P. H. Rabinowitz,
Minimax Methods in Critical Point Theory With Applications to Differential Equations,
CBMS Regional Conference Series in Mathematics,  Amer Mathematical Society, 1986.

\bibitem[S1]
M. Shibata, 
Stable standing waves of nonlinear Schr\"odinger equations with a general nonlinear term,
Manuscripta Math. 143 (2014), no. 1-2, 221--237.

\bibitem[S2]
M. Shibata, 
A new rearrangement inequality and its application for $L^2$-constraint minimizing problems,
Math. Z. 287 (2017), no. 1-2, 341--359.

\bye


\vfil
\break

{

\baselineskip 30 pt

\fontscale[1100]

\def\cmd#1#2{\line{\indent\hbox to 2.5 cm{\tt \string#1\hfil}\hskip 1 cm $\displaystyle #1$\hfil}}

\let\def=\cmd

\centerline{\bf --- macros ---}

\bigskip

\input macro.mac

}


\bye